\documentclass[preprint,12pt]{elsarticle}

\usepackage{amssymb,amsmath,verbatim,xcolor}
\usepackage{todonotes}
\usepackage{bm}
\usepackage{url}

\newdefinition{remark}{Remark}

\usepackage{float}
\def\bn{{\bf n}}

\def\ljump{{[\![}}
\def\rjump{{]\!]}}
\def\lavg{{\{\hspace{-.05in}\{}}
\def\ravg{{\}\hspace{-.05in}\}}}

\newcounter{bla}

\journal{Computer Physics Communications}

\begin{document}
\begin{frontmatter}
\title{An Efficient High-order Numerical Solver for Diffusion Equations with Strong Anisotropy}
\tnotetext[t1]{The first and third authors are supported by the Office of Science of the U.S. Department of Energy under Contract No. DE-AC05-00OR22725..}

\author{David Green\fnref{label1}}
\ead{greendl1@ornl.gov}
\author{Xiaozhe Hu\fnref{label2}}
\ead{Xiaozhe.Hu@tufts.edu}
\author{Jeremy Lore\fnref{label1}}
\ead{lorejd@ornl.gov}
\author{Lin Mu\corref{cor1}\fnref{label3}}
\ead{linmu@uga.edu}
\author{Mark L. Stowell \fnref{label4}}
\ead{stowell1@llnl.gov}
\cortext[cor1]{Corresponding author}
\fntext[label1]{Fusion \& Materials for Nuclear Systems Division
Oak Ridge National Laboratory, Oak Ridge, TN 37831, USA}
\fntext[label2]{Department of Mathematics, Tufts University, Medford, MA 02155}
\fntext[label3]{Department of Mathematics, University of Georgia, Athens, GA 30602}
\fntext[label4]{Center for Applied Scientific Computing at Lawrence Livermore National Laboratory, 7000 East Ave., Livermore, CA 94550}

\date{}

\begin{abstract}
In this paper, we present an interior penalty discontinuous Galerkin finite element scheme for solving diffusion problems with strong anisotropy arising in magnetized plasmas for fusion applications. We demonstrate the accuracy produced by the high-order scheme and develop an efficient preconditioning technique to solve the corresponding linear system, which is robust to the mesh size and anisotropy of the problem. Several numerical tests are provided to validate the accuracy and efficiency of the proposed algorithm.
\end{abstract}
\begin{keyword}
Anisotropic diffusion equation, Interior penalty discontinuous Galerkin, High-order method, Subspace correction methods, Iterative methods
\end{keyword}

\end{frontmatter}

\section{Introduction}
Anisotropic diffusion is a common physical phenomenon and describes processes where the diffusion of some scalar quantity is direction dependent. 
In fusion plasmas diffusion tensors can be extremely anisotropic due to the high temperature and large magnetic field strength. 
For example, in plasma magnetic fusion devices, the electron heat conductivity along the magnetic field can be greater than that across the magnetic field by a factor of $10^6$ or even may reach the order of $10^{12}$ \cite{vEKdB2014}. This anisotropy is due to the fact that the gryomotion of charged particles in a magnetic field results in slow transport perpendicular to the field while particles can travel comparatively long distances parallel to the field before undergoing a collision.

Strongly anisotropic diffusion problems present a numerical challenge, because  errors in the direction parallel to the magnetic field may have significant effect on transport in the perpendicular direction. The widely used strategy to solve the anisotropic diffusion equations is to use finite-difference schemes with meshes aligned with the magnetic field direction. Flux aligned coordinates are used in the fusion community to obtain accurate diffusion simulations \cite{DUXSW2009}. However, the generation of aligned meshes is still limited to comparably simple magnetic field topologies. For practical geometries, aligned mesh generation is still an open problem, particularly for fusion confinement devices with complicated magnetic topologies and plasma-facing component geometries. Alternative discretization schemes which do not rely on a magnetic field aligned grid have also been studied \cite{O2011,vEKdB2014}, and such a discretization is the focus of this paper. 

In order to relax the difficulties in generating the aligned mesh, there is increasing interest in non-aligned mesh schemes. These methods must overcome recognized difficulties arising for strongly anisotropic diffusion problems with non-aligned meshes, which include numerical perpendicular pollution, non-positivity, and loss of convergence \cite{vEKdB2014}. Finite difference methods adopting interpolations aligned to the parallel diffusion direction have been proposed in \cite{SSGLTS2020}. Proper flux construction and high order finite difference schemes on the non-aligned meshes have been proposed in \cite{GYKL2005,GLT2007}. Due to the flexibility in the geometry, the finite volume methods on the non-aligned meshes have also been investigated for example in~\cite{HFS2013} and the high order finite volume schemes has been studied in \cite{CKL2015}. In the finite element community, the continuous linear finite element coupled with M-matrix adaption has been studied in \cite{LH2010}. Jardin \cite{J2004} applies a finite element method with reduced quintic triangular finite elements where the quintic basis functions are constrained to enforce $C^1$ continuity across element boundaries. 

Compared with other finite element methods, discontinuous Galerkin methods have many advantages, including flexibility in mesh generation, adopting hp-adaptive refinement, and preserving the essential conservation properties. Such methods have also been used to discretize the anisotropic diffusion operator by a discontinuous Galerkin scheme in \cite{HWS2016} and a hybrid discontinuous Galerkin scheme in \cite{GBSST2020} for the non-aligned meshes. Both of these two approaches are able to achieve high order approximation by increasing the polynomial order which  usually results in a similarly accurate solution with fewer total degrees of freedom than lower-order elements. In this paper, we shall also employ a high order scheme and additionally demonstrate the advantages in the interior penalty discontinuous Galerkin method. More details on the different types of discontinuous Galerkin methods can be found in \cite{ABCM2002}.

However, the discontinuous Galerkin method will result in even a larger number of degrees of freedom than continuous finite element method and the computational cost can be higher due to the resulting ill-conditioned linear system. In fact, the condition number is even worse in the strong anisotropic diffusion case. As a result, an efficient implementation requires the use of advanced iterative methods for solving the resulting linear systems. Efficient solvers arising from the discontinuous Galerkin method can be found in \cite{AA2007} (Schwarz method), \cite{DLVZ2006} (multilevel method) and \cite{BZ2005} (multigrid methods). The subspace correction method \cite{X1992} and auxiliary space preconditioning method \cite{X1996} have also been used to construct an efficient linear solvers for the discontinuous Galerkin methods, see \cite{DZ2009,ASVZ2017}. In this paper, we will use the auxiliary space method to construct a robust preconditioner.

Our particular application focus for this paper is magnetized plasmas for fusion applications, where the scale separation induced by the magnetic field produces conduction coefficients that are several orders of magnitude larger in the parallel direction (denoted as $D_{\|}$) than in the perpendicular direction (denoted as $D_{\perp}$). This work demonstrates that a high-order method is able to provide satisfactory numerical solutions on non-aligned meshes. We shall also construct an effective and robust preconditioner to solve the corresponding linear system.
The remainder of this paper is organized as follows. Section~\ref{sect:FEM} reviews the notation and proposes the numerical schemes for the anisotropic diffusion equations. The numerical experiments are reported in Section~\ref{sect:num-accuracy} to validate the accuracy tests for our proposed schemes. Then Section~\ref{sect:preconditioner} and Section~\ref{sect:num-preconditioner} are contributed to develop and verify the efficient linear solver and numerical tests for problems with strong anisotropy. The concluding remarks and future research plans will be discussed in Section~\ref{sect:conclusion}.

\section{Finite Element Spaces and Finite Element Schemes}\label{sect:FEM}
In this section, we will review the notations in the continuous Galerkin (CG) and discontinuous Galerkin (DG) finite element scheme in the later sections. The DG methods will be used to discretize the given equations and the CG methods will be used in constructing the preconditioner.

\subsection{Preliminaries}
The problem under consideration is the following two-dimensional equation with anisotropic diffusivity. Let $\Omega$ be a bounded connected domain in $\mathcal{R}^2$ with Lipschitz boundary $\partial\Omega$ and let $f\in L^2(\Omega)$.
We consider the following steady state anisotropic diffusion equation:
\begin{eqnarray}
-\nabla\cdot(\mathbb{D}\nabla u) &=& f, \mbox{ in }\Omega, \label{eqn: model-1}\\
u &=& 0,\text{ on }\partial\Omega, \label{eqn:model-2}
\end{eqnarray}
where 
the diffusion coefficient tensor is given by 
\begin{eqnarray}
\mathbb{D} = 
\begin{pmatrix}
b_1 &-b_2\\b_2 &b_1
\end{pmatrix}
\begin{pmatrix}
D_{\|} &0\\ 0 & D_{\perp}
\end{pmatrix}
\begin{pmatrix}
b_1 &b_2\\-b_2 &b_1
\end{pmatrix}.
\label{eq:diff-coef}
\end{eqnarray}
The direction of the anisotropy, or the magnetic field, is given by a unit vector ${\bf b} = (b_1,b_2)^\top$. Here $D_{\|}$ and $D_{\perp}$ represent the parallel and the perpendicular diffusion coefficient. Moreover, the anisotropy level is such that $D_{\|}$ is several orders of magnitude larger than $D_{\perp}.$

Throughout this paper, the standard notations for Sobolev spaces will be adopted. For example, for a bounded domain $\omega\subset\mathcal{R}^2$, we denote by $H^m(\omega)$ the standard Sobolev space of order $m\ge 0$ with norm and semi-norms $\|\cdot\|_{m,\omega}$ and $|\cdot|_{m,\omega}$ respectively. If $m = 0$ and/or $\omega=\Omega$, we shall omit $m$ and/or $\omega$ in the norm notations.

\vskip.1in
\noindent\textbf{Partition}: Denote a shape-regular partition of domain $\Omega$ into two-dimensional simplices $T$ (triangles or quadrilaterals) as  $\mathcal{T}_h:=\cup T$. 
Denote $\mathcal{E}_h=\cup_{T\in\mathcal{T}_h} \partial T$ as the interfaces for all the elements $T$ with $\mathcal{E}_h^0$ denoting the interior edges. Let $h = \max_{T\in\mathcal{T}_h} (h_T)$, where $h_T$ denotes the diameter of element $T$.

\vskip.1in
\noindent\textbf{Trace Operators}:
For piecewise functions, we further introduce the jumps and averages as follows. For any edge $e=T^+\cap T^-\in\mathcal{E}_h^0$, with $\bn^\pm$ as the outward unit normal to $\partial T^\pm$, the jumps of a scalar valued function $u$ across $e$ are defined as
\begin{eqnarray*}
\ljump u \rjump = u^- - u^+,\
\end{eqnarray*}
and the averages of a scalar valued function $u$ 
is defined as
\begin{eqnarray*}
\lavg u \ravg=\frac{1}{2}\left(u^++u^-\right).  
\end{eqnarray*}
On the boundary $\partial\Omega$, when $e\in\partial\Omega$,  the jump and average operators are defined as 
$\lavg u \ravg = u|_e$, 
$\ljump u \rjump = u|_e$.

\vskip.1in
\noindent\textbf{Continuous Finite Element Space}: Let 
\begin{eqnarray}
V_\text{CG} = \{v\in H_0^1(\Omega)|\; v|_T\in P_k(T),\forall T\in\mathcal{T}_h\},
\end{eqnarray}
where $P_k(T)$, $k\ge 1$, denotes the polynomials with degree less than or equal to $k$.

\vskip.1in
\noindent\textbf{Discontinuous Finite Element Space}: Let 
\begin{eqnarray}
V_\text{DG} = \{v\in L^2(\Omega)|\; v|_T\in P_k(T),\forall T\in\mathcal{T}_h\}.
\end{eqnarray}


\subsection{Continuous Galerkin Finite Element Method}
In this section, we shall introduce the $H^1$-finite element method (FEM). 
The $H^1$-FEM is to find the numerical solution $u_h\in V_\text{CG}$, such that
\begin{eqnarray}
A_\text{CG}(u_h,v) = \sum_{T\in\mathcal{T}_h}\int_T fv \, \mathrm{d}T,\ \forall v\in V_\text{CG},
\end{eqnarray}
where the bilinear form
\begin{eqnarray}
A_\text{CG}(w,v) = \sum_{T\in\mathcal{T}_h}\int_T\mathbb{D}\nabla w\cdot\nabla v \, \mathrm{d}T.
\end{eqnarray}

\begin{remark}
Well-posedness of CG for isotropic diffusion equations can be found in the many literatures, for example \cite{BS2008}. It is noted that the CG scheme can produce a symmetric positive definite linear system.
\end{remark}

\begin{remark}
As with the traditional analysis for the isotropic diffusion equation, we can expect the optimal convergence rate for $L^2$-error at the order $\mathcal{O}(h^{k+1})$ if $u\in H^{k+1}(\Omega)$. However, in the strong anisotropic case, the optimal rate of convergence is usually hard to obtain with low order schemes due to spurious numerical diffusion. Some benchmark tests can be found in \cite{HH2008}. 


\end{remark}

\begin{remark}The following condition number estimation result for two dimensional anisotropic diffusion equations discretized by the linear CG method can be found in \cite{KHX2014}.
The condition number of the stiffness matrix for $k = 1$ is bounded by
\begin{eqnarray*}
\kappa(A)\le CN\times\bigg(\frac{1}{\lambda_{\min}}\max_j\sum_{T\in\omega_j}|T|\|(F_T')^{-1}\mathbb{D}_T(F_T')^{-\top}\|_2 \bigg)
\times
\bigg(1+\ln\frac{|\bar{T}|}{|T_{\min}|}\bigg).
\end{eqnarray*}
Here $C>0$ denotes a constant, $\lambda_{\min}$ denotes the minimum eigenvalue of $\mathbb{D}$ on $\Omega$, $N$ denotes the number of elements, $\omega_j$ is the element patch associated with the $j$th vertex, $\mathbb{D}_T$ is the average of $\mathbb{D}$ over $T$, $F_T'$ is the Jacobian matrix of the affine mapping from the reference element to the physical element $T$, and $|\bar{T}| = \frac{1}{N}|\Omega|$ denotes the average element size. In this conditioning estimate, the first factor $N$ corresponds to the condition number of the stiffness matrix for the Laplacian operator. The second factor reflects the effects of the volume weighted, combined alignment and equidistribution quality measure of the mesh with respect to anisotropic diffusion tensor $\mathbb{D}$. The third term measures the mesh volume-nonuniformity. We also note that there are many references (e.g., \cite{LH2010,H2011,LH2013}) developing the anisotropy-adaptive ($\mathbb{D}^{-1}$-adaptive) meshing techniques in the continuous linear finite element framework. However, it is not a focus for this paper.
\end{remark}

\subsection{Discontinuous Galerkin Finite Element Method}

Now, we are ready to introduce the numerical algorithm for solving anisotropic diffusion equation by using piecewise polynomials without any continuous constraints.
The interior penalty discontinuous Galerkin (IPDG) numerical algorithm is to find $u_h\in V_\text{DG}$ such that
\begin{eqnarray}\label{eqn:DG}
A_\text{DG}(u_h,v) = 
&&\sum_{T\in\mathcal{T}_h}\int_T fv \, \mathrm{d}T,\; \forall v\in V_\text{DG},
\end{eqnarray}
where
\begin{eqnarray*}
A_\text{DG}(u_h,v) =&& \sum_{T\in\mathcal{T}_h}\int_T\mathbb{D}\nabla u_h\cdot\nabla v \, \mathrm{d}T -\sum_{e\in\mathcal{E}_h}\int_e\lavg\mathbb{D}\nabla u_h\cdot\bn\ravg\ljump v \rjump \, \mathrm{d}s \\
&-& \beta\sum_{e\in\mathcal{E}_h}\int_e\lavg\mathbb{D}\nabla v\cdot\bn\ravg\ljump u_h \rjump \, \mathrm{d}s 
+ \sum_{e\in\mathcal{E}_h}\int_e\alpha \ljump u_h\rjump\ljump v \rjump \, \mathrm{d}s.
\end{eqnarray*}
If $\beta = 1$, we obtain the symmetric IPDG scheme and $\beta = -1$, we obtain the non-symmetric IPDG scheme. The parameter $\alpha$ has to be chosen to ensure stability.

\begin{remark}
The well-posedness of DG can be found in \cite{SR2011}. Briefly speaking, in the case of symmetric IPDG scheme, the penalty parameter $\alpha$ is required to be large enough to ensure stability; however, in the case of non-symmetric IPDG scheme, we only need the penalty parameter to be positive.  
\end{remark}

\begin{remark}
For two-dimensional simulation, the penalty parameter $\alpha$ for symmetric IPDG schemes with $\mathbb{D} = \mathbb{I}$  is widely chosen as
\begin{eqnarray*}
\alpha = \frac{4k(k+1)}{h}.
\end{eqnarray*}
In the remainder of this paper, we shall only focus on the symmetric IPDG scheme with the penalty parameter $\alpha$ below,
\begin{eqnarray}
\alpha_1 = \frac{4k(k+1)D_{\|}}{h},\ \alpha_2 = \frac{4k(k+1)}{h}\bn\cdot(\mathbb{D}\bn).\label{eq:alpha_val}
\end{eqnarray}
In the above parameter setting, $\alpha_1$ is weighted with the nodal distance in each element and $\alpha_2$ is based on the modulation of the stabilization with respect to the actual strength of the parallel diffusion in the normal direction of each face.
\end{remark}


\begin{remark}
Let $\lambda_{\max} :=\max(1,\lambda_{\max,\mathbb{D}})$, where $\lambda_{\max,\mathbb{D}}$ indicates the maximum eigenvalue of $\mathbb{D}$ on $\Omega$. Let $\lambda_{\min}$ denotes the minimum eigenvalue of $\mathbb{D}$ on $\Omega$. Suppose $u\in H^{k+1}(\Omega)$ and the numerical scheme is wellposed, then the following error estimate \cite{ESZ2009} holds,
\begin{eqnarray}
\| u-u_h\|&\le& C \frac{\lambda_{\max}}{\lambda_{\min}}h^{k+1}\|u\|_{k+1}
\end{eqnarray}
The above inequality also demonstrates that in the case of strong anisotropy, there may be large numerical pollution in the low order schemes due to the $\lambda_{\max}$ term.
\end{remark}

\section{Numerical Experiments for Testing Accuracy}\label{sect:num-accuracy}
In this section, we first study the accuracy of the above IPDG numerical schemes. We demonstrating that:
\begin{itemize}
\item  The simulation on the aligned meshes can resolve the anisotropy to some level and is able to provide the satisfactory numerical solutions.
\item The high-order IPDG schemes on the un-aligned meshes  are still able to produce reliable numerical simulations. 
\end{itemize}

Errors are measured in $L^2$-norm to validate our results, with the optimal rates in convergence being at the order $\mathcal{O}(h^{k+1})$, where $k$ denotes the polynomial degree in the IPDG scheme.

\subsection{Aligned Mesh}
In this section we demonstrate the IPDG scheme on a series of aligned meshes of increasing resolution for three problems of increasing solution variability in the perpendicular direction. Section~\ref{Num-Test2} examines the scheme's performance for the non-aligned meshes. 

\subsubsection{Test: Error versus $D_{\parallel}$ for varying mesh and penalty parameters}\label{Num-Test1}
\begin{figure}[H]
\centering
\begin{tabular}{ccc}
\includegraphics[width =.3\textwidth]{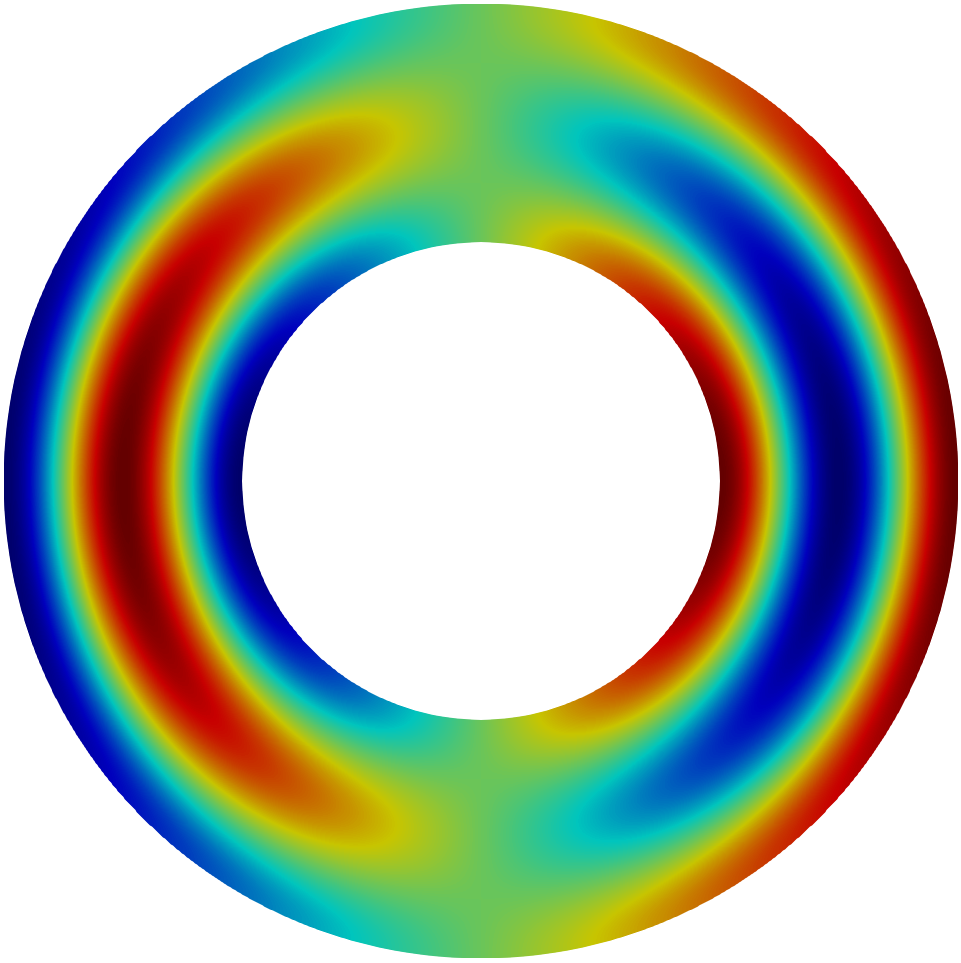}
&\includegraphics[width =.3\textwidth]{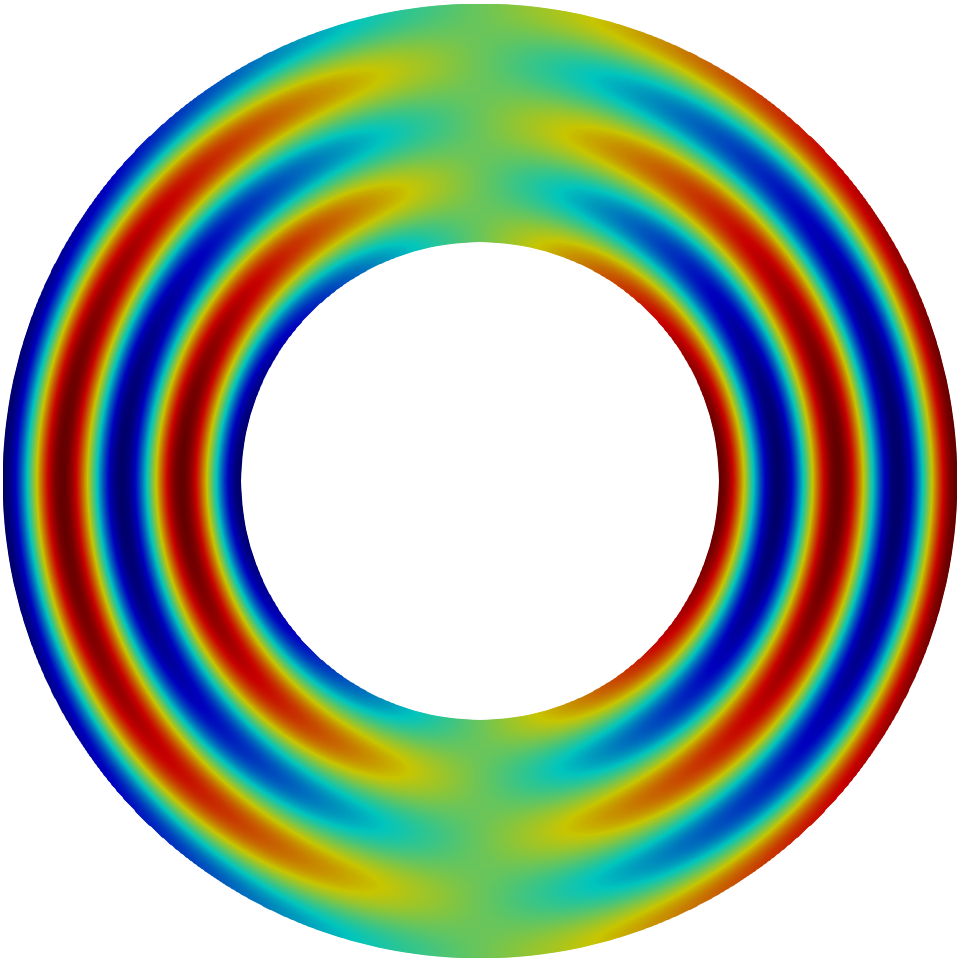}
&\includegraphics[width =.3\textwidth]{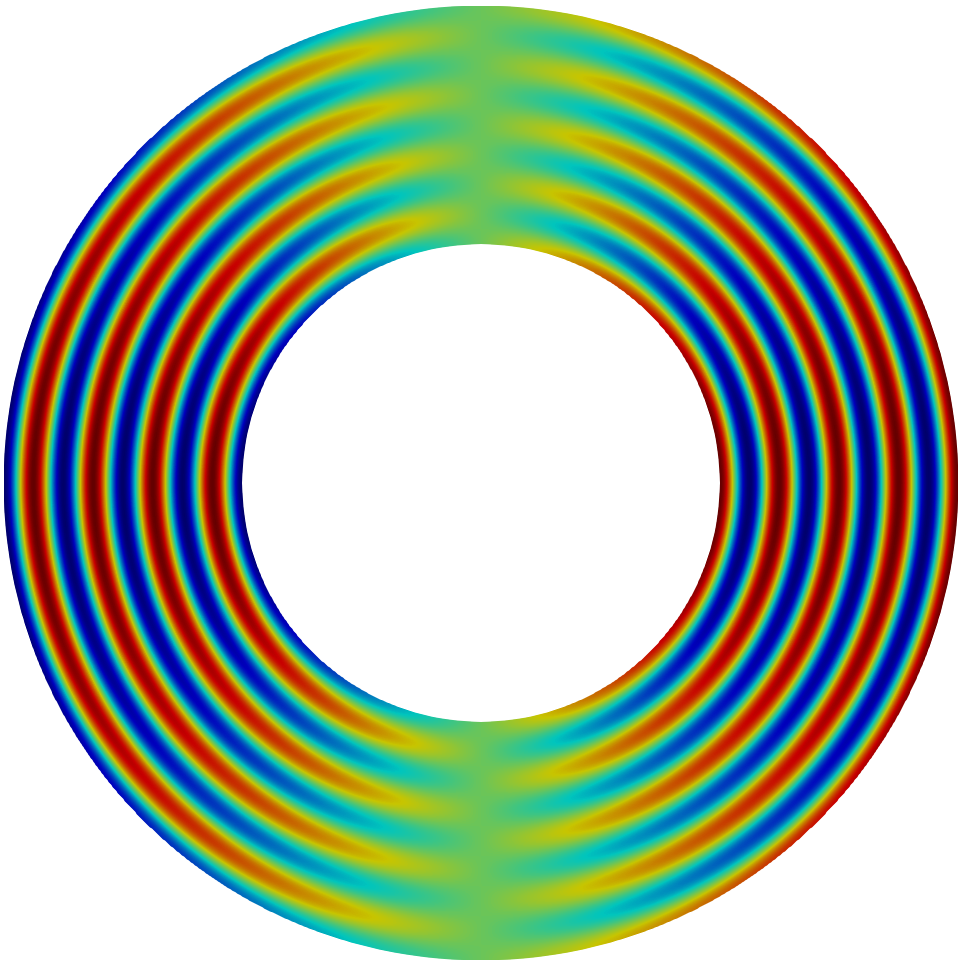}\\
(a) & (b) &(c)
\end{tabular}
\caption{Test~\ref{Num-Test1}: Analytical solutions with respect to different radial (perpendicular) oscillation frequencies: (a) $\omega = 1$; (b) $\omega = 2$; (c) $\omega = 4$.}\label{fig:annulus-1}
\end{figure}
\noindent\textbf{Problem Setting}: In this test, the computational domain $\Omega$ is an annulus with exterior radius $r = 2$ and interior radius $r = 1$, and the magnetic field direction is circular with the following expression
\begin{eqnarray}
u &=& \cos(2\pi\omega r)\cos(\theta),\notag \\
b_1 &=& \dfrac{y}{r}, \ b_2 = -\dfrac{x}{r}\notag\\
f &=& \frac{4x}{r^3}\left( (\pi^2r^2\omega^2 +  \frac{D_{\|}}{4})\cos(2\pi\omega r) + \frac{\pi\omega r}{2}\sin(2\pi\omega r)\right).\label{ex:TestAnnulus}
\end{eqnarray}
The diffusion tensor (\ref{eq:diff-coef}) is set
with 
$D_{\perp} = 1.0$ and varying values of $D_{\|}$. Here, $\omega$ is the frequency of oscillation in the perpendicular direction.

\begin{figure}[H]
\centering
\begin{tabular}{ccc}
\includegraphics[width =.3\textwidth]{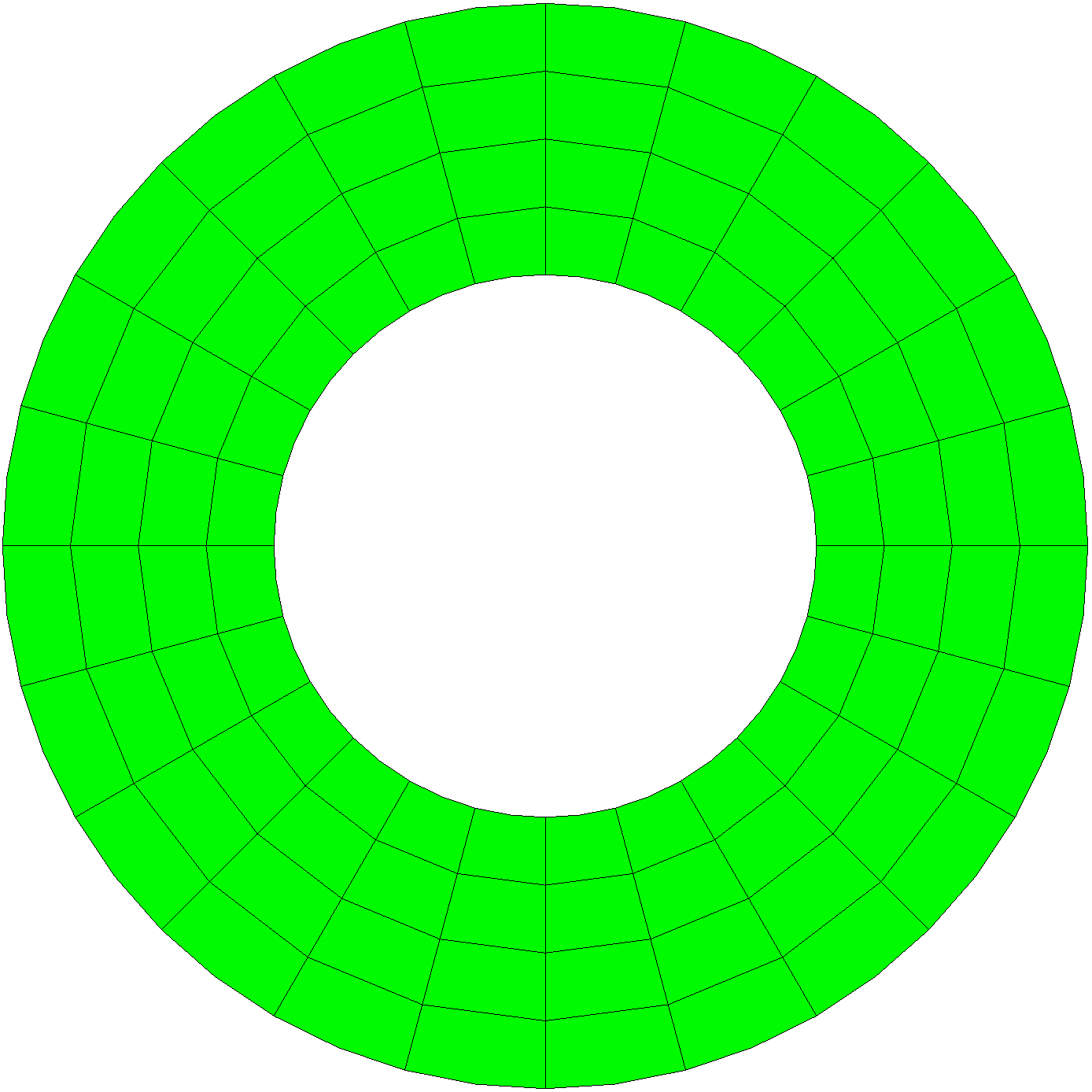}
&\includegraphics[width =.3\textwidth]{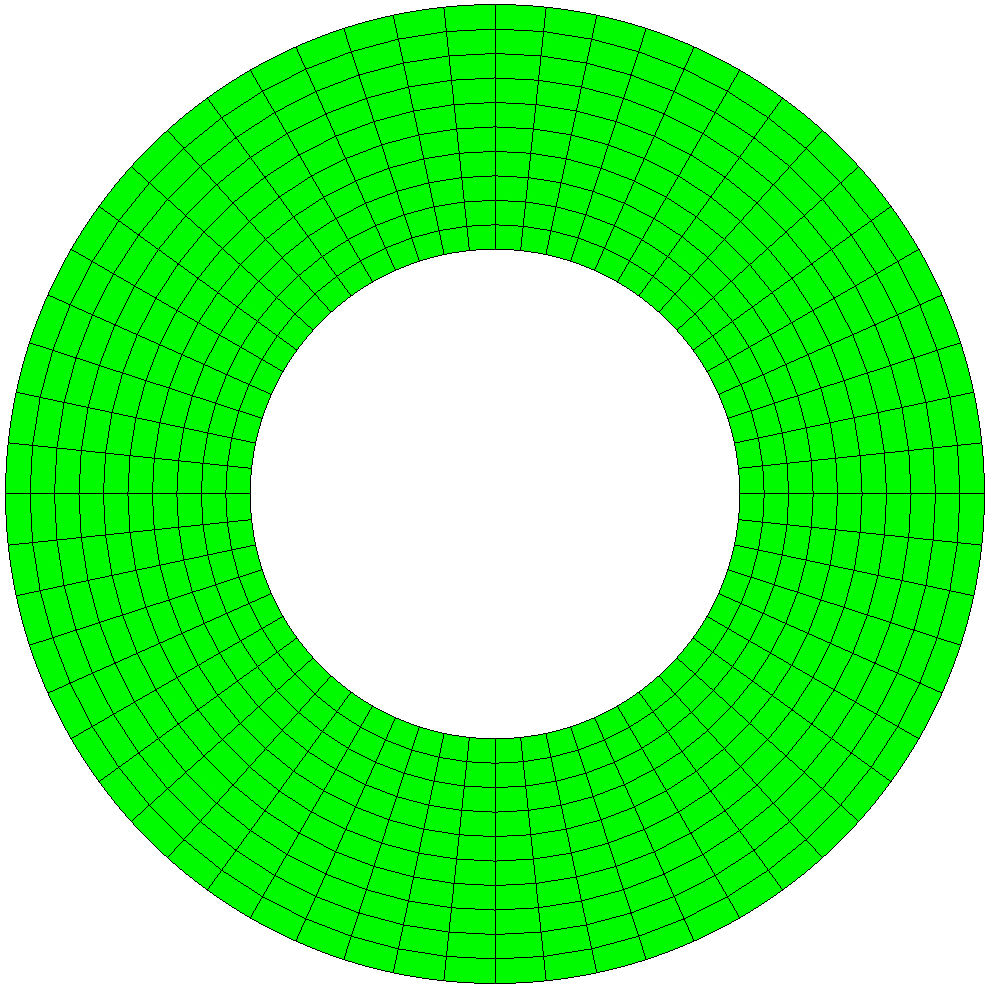} 
&\includegraphics[width =.3\textwidth]{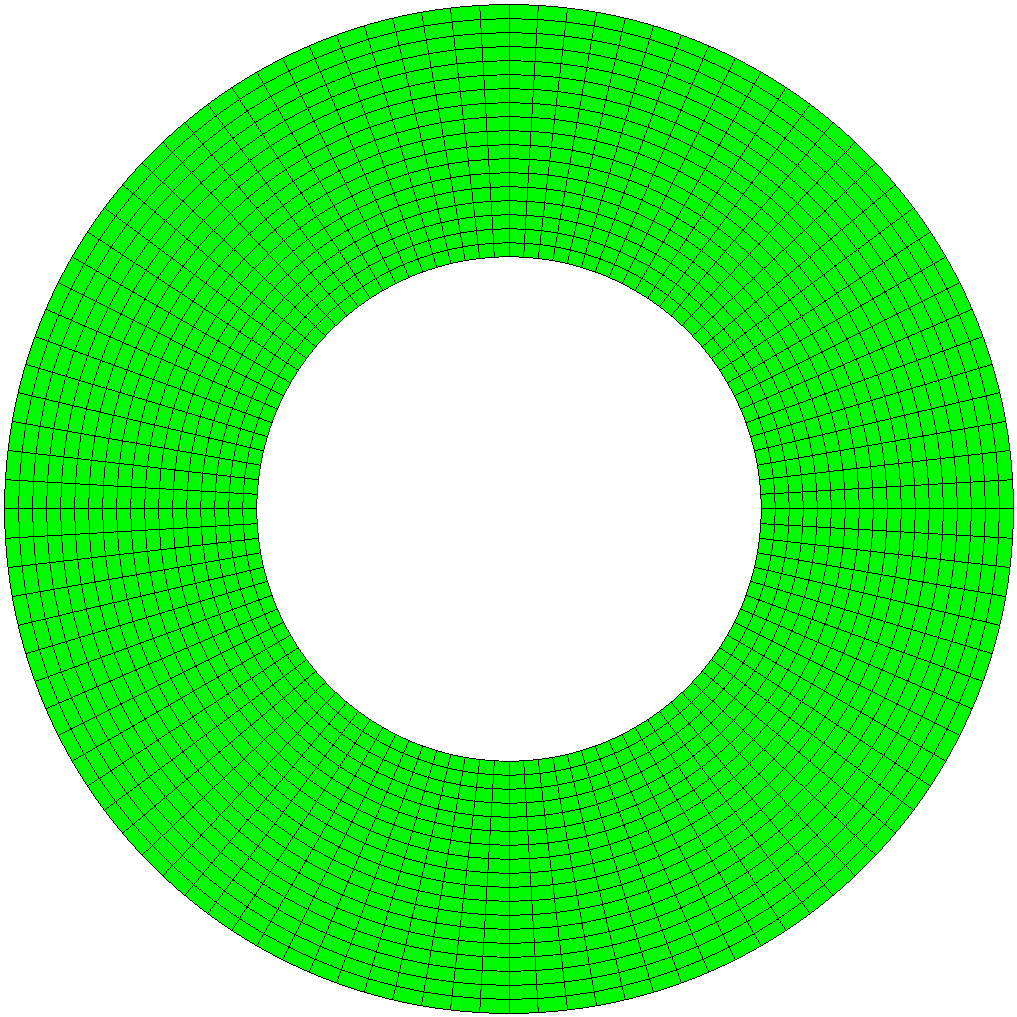}\\
(a) & (b) & (c)
\end{tabular}
\caption{Test~\ref{Num-Test1}: Illustration of computational meshes (a) $N_r = 4$; (b) $N_r = 10$; (c) $N_r=18$.}\label{fig:annulus-2}
\end{figure}

The manufactured solutions are plotted in Figure~\ref{fig:annulus-1} for $\omega = 1,2,4$ and the discretizations are based on the quadrilateral meshes shown in Figure~\ref{fig:annulus-2}. For $\omega = 1$, we utilize the mesh with $N_r = 4$ (Figure~\ref{fig:annulus-2}a); for $\omega = 2$, $N_r = 10$ (Figure~\ref{fig:annulus-2}b); and for $\omega = 4$, $N_r = 18$ (Figure~\ref{fig:annulus-2}c).

\begin{figure}[H]
\centering
\begin{tabular}{cc}
\includegraphics[width =.45\textwidth]{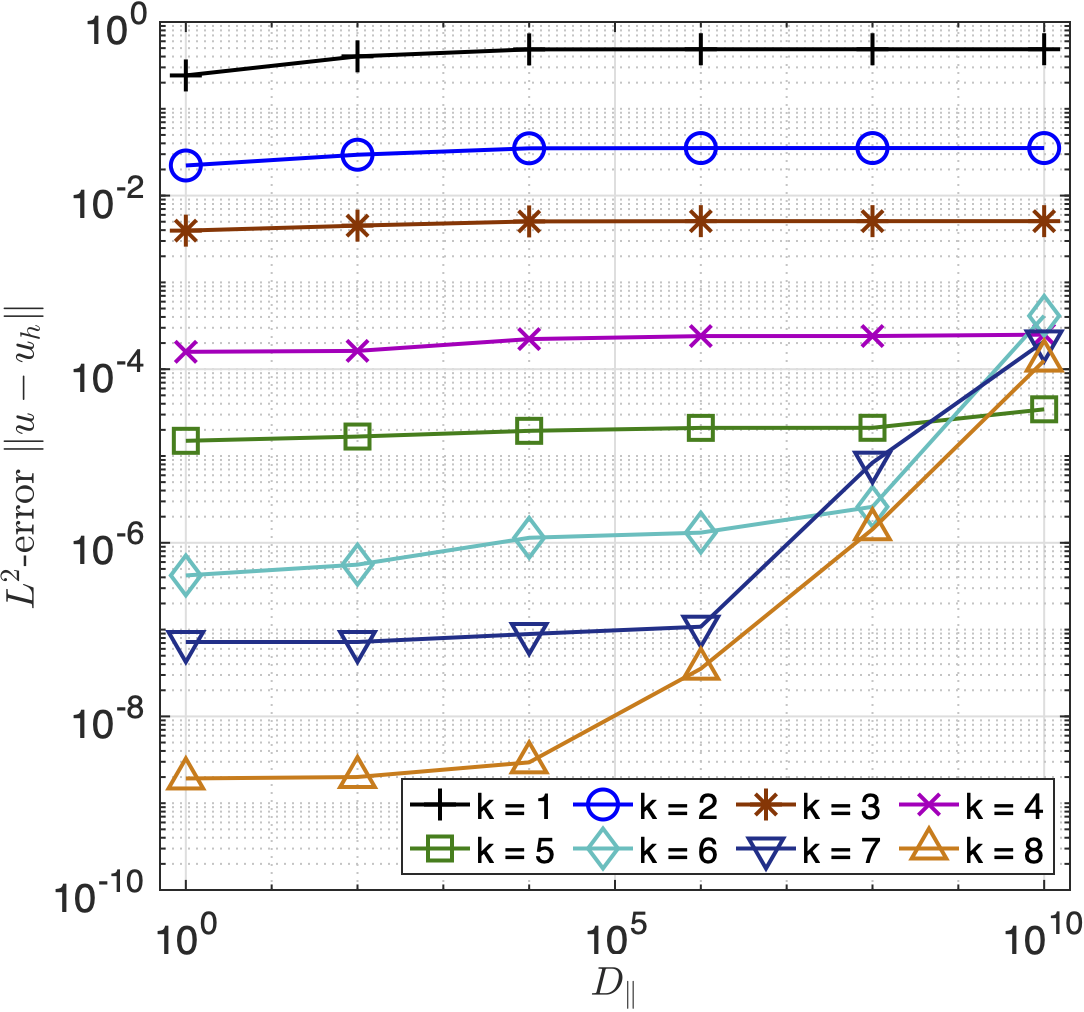}
&\includegraphics[width =.45\textwidth]{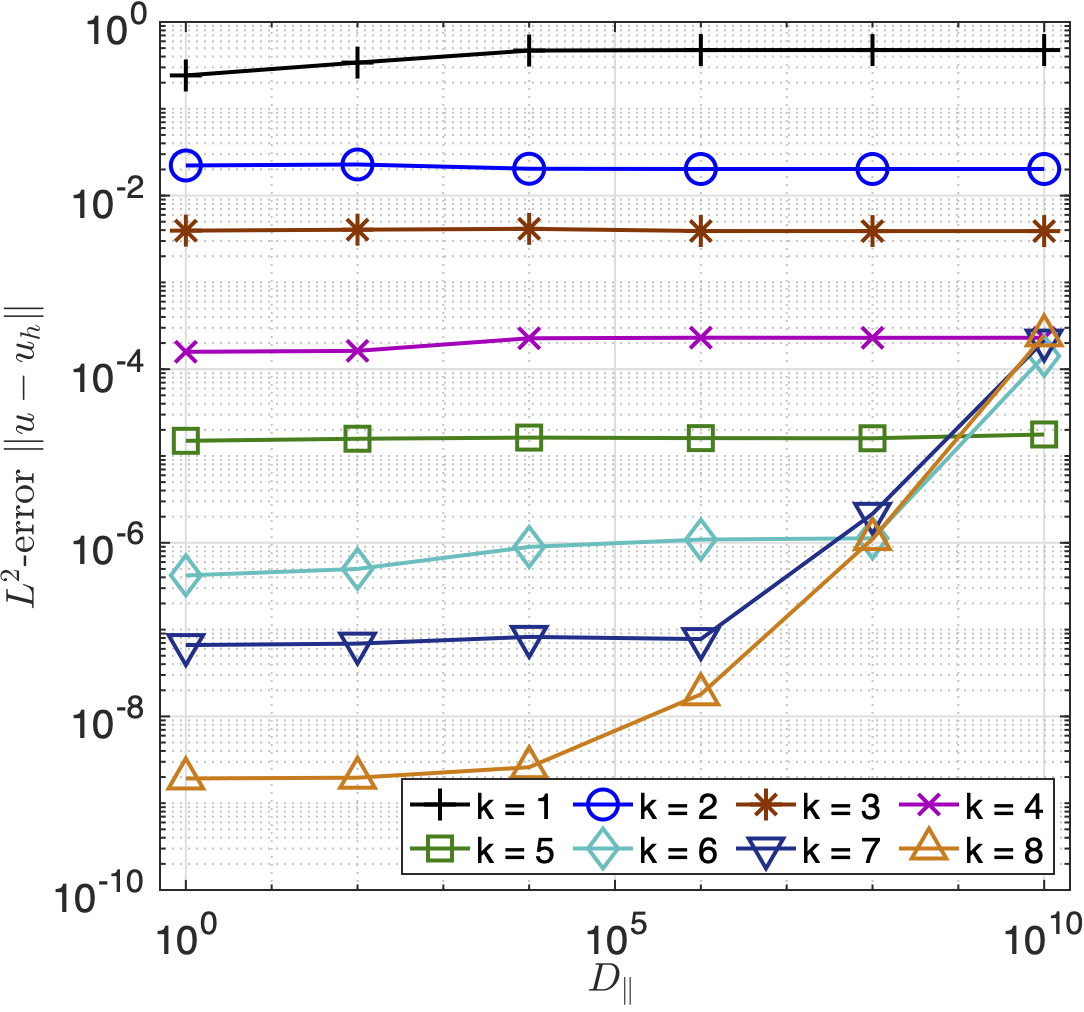}\\
(a) & (b) 
\end{tabular}
\caption{Test~\ref{Num-Test1}: Illustration of $L^2$-error as a function of the anisotropy ratio $D_{\|}$ for $\omega = 1$ and $k = 1,\cdots,8$ on the mesh in Figure~\ref{fig:annulus-2} (a) with different penalty parameters: (a) $\alpha_1$; (b) $\alpha_2$.}\label{fig:annulus-3}
\end{figure}

The numerical performance for $\omega = 1$ on the Mesh with $N_r = 4$ (Figure~\ref{fig:annulus-2}a) is summarized in Figure~\ref{fig:annulus-3}. The corresponding linear system is solved directly with package ``sparsecholesky" \cite{eigenweb}. The 1D plots show the values of the $L^2$-error for $D_{\|}$ varying between $10^0$ and $10^{10}$. With both penalty parameters $\alpha_1$ and $\alpha_2$, the value of the $L^2$-error remains constant whatever the value of $D_{\|}$, except for $k = 6,7,8$ where it finishes by increasing at large $D_{\|}$ values ($D_{\|}\ge 10^5$). This behavior has to be related to the constant increase of the condition number of the matrix with $D_{\|}$ when solving the linear system with a direct solver. 
As long as the anisotropy remains moderate in our test case, the $L^2$-error remains dominated by the polynomial interpolation errors. However, in the case with strong anisotropy, the error from the linear solver dominates over the interpolation error, and thus leads to an increase of the $L^2$-error. The error behavior in the high-order scheme may be solved by employing other direct solvers.  For 
either value of $\alpha_i$ ($i=1,2$), the behavior is very similar when comparing the $L^2$-error.
For simplicity, we shall choose the $\alpha_1$ penalty parameter in all of the following tests.

\begin{figure}[H]
\centering
\begin{tabular}{cc}
\includegraphics[width =.45\textwidth]{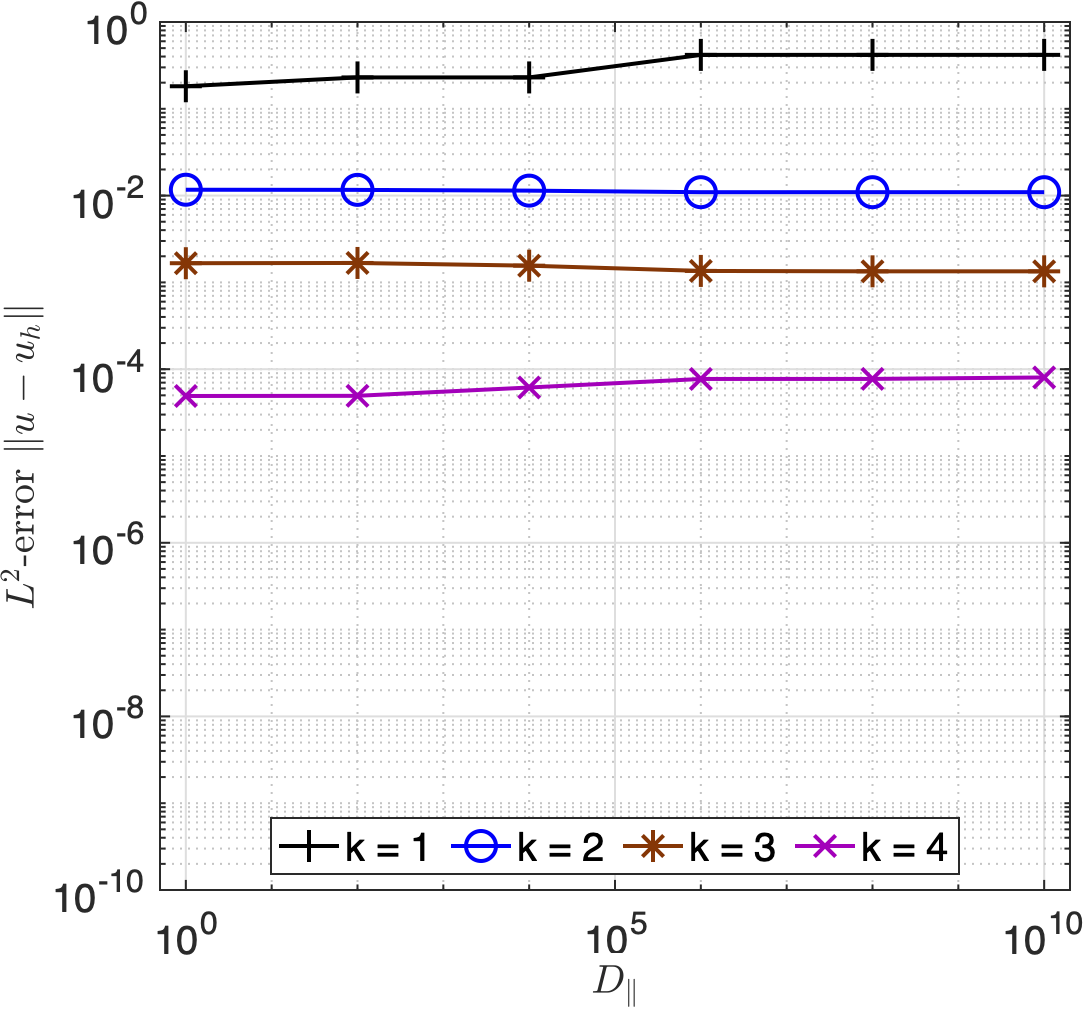}
&\includegraphics[width =.45\textwidth]{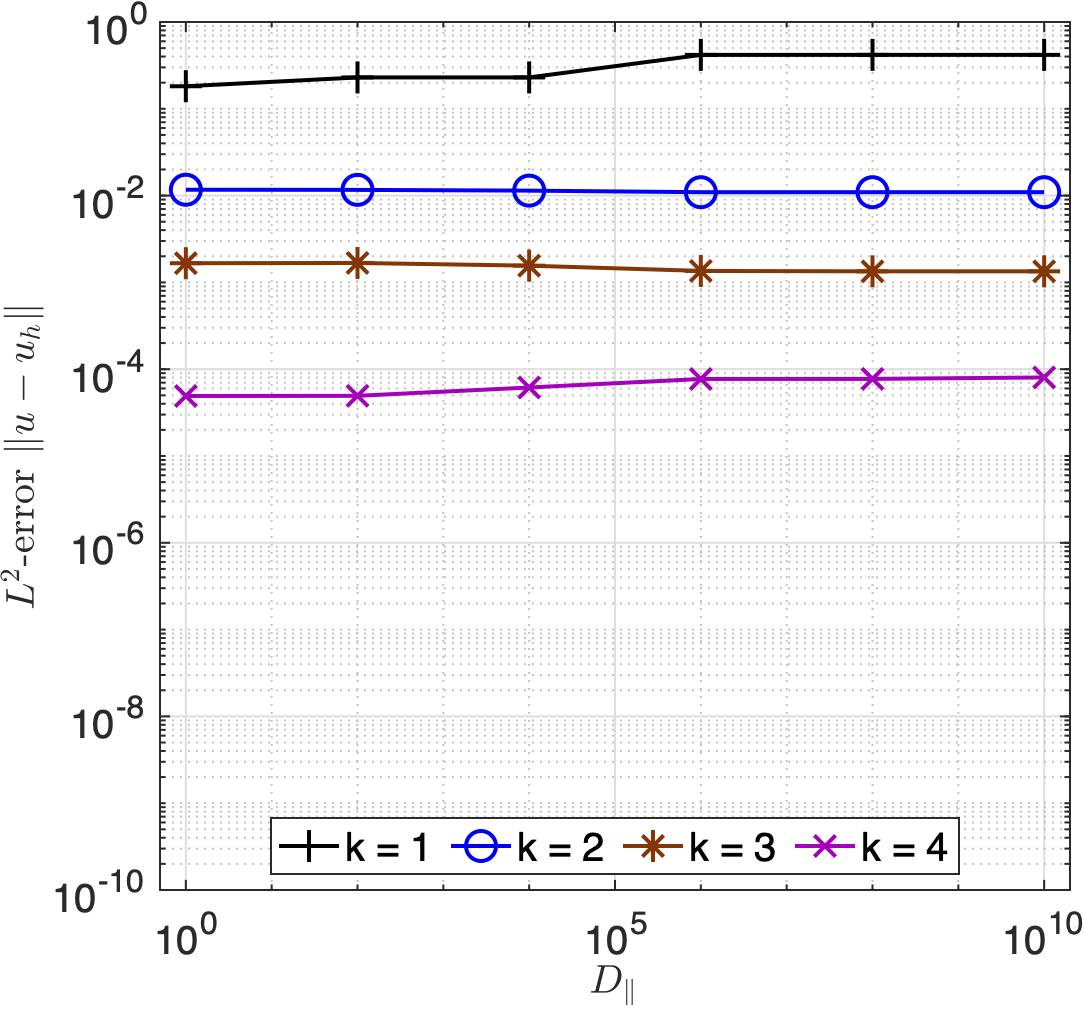}\\
(a) & (b) 
\end{tabular}
\caption{Test~\ref{Num-Test1}: Illustration of $L^2$-error as a function of the anisotropy ratio $D_{\|}$ for $k = 1,\cdots,4$ with $\alpha_2$ on: (a) Mesh in Figure~\ref{fig:annulus-2}(b) with $\omega = 2$; (b) Mesh in Figure~\ref{fig:annulus-2}(c) with $\omega = 4$.}\label{fig:annulus-4}
\end{figure}

Next, we shall test the performance for $\omega = 2$ and $\omega = 4$ with $ k = 1,2,3,4$. Since our direct solver can handle the condition numbers in these linear systems, we can expect flat behavior in the error with $D_{\parallel}$. Because of the increase in perpendicular frequency we partition the domain with more divisions in the radial direction as shown in Figure~\ref{fig:annulus-2}a-b. The $L^2$-errors with respect to the values in $D_{\|}$ are plotted in Figure~\ref{fig:annulus-4}. Again, we obtain the expected performance that the $L^2$-errors remain almost constant with $D_{\|}$ and thus validate our stabilization conclusions in the parameter $\alpha_1$.

\subsection{Non-Aligned Mesh}\label{Num-Test2}
In this section, we shall demonstrate the performance of the high-order scheme on the non-aligned mesh and show that satisfactory numerical solutions can be obtained.
\subsubsection{Test: Constant magnetic field}\label{Num-Test2-1}
\begin{figure}[H]
\centering
\begin{tabular}{cc}
\includegraphics[width =.4\textwidth]{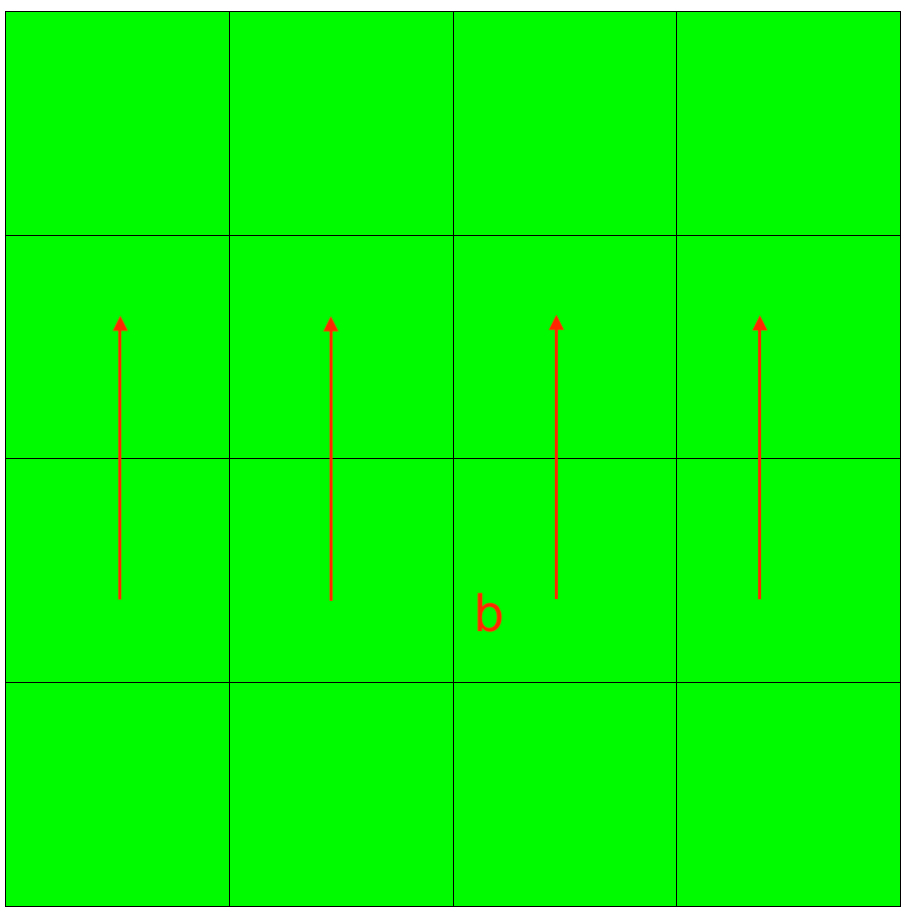}
&\includegraphics[width =.4\textwidth]{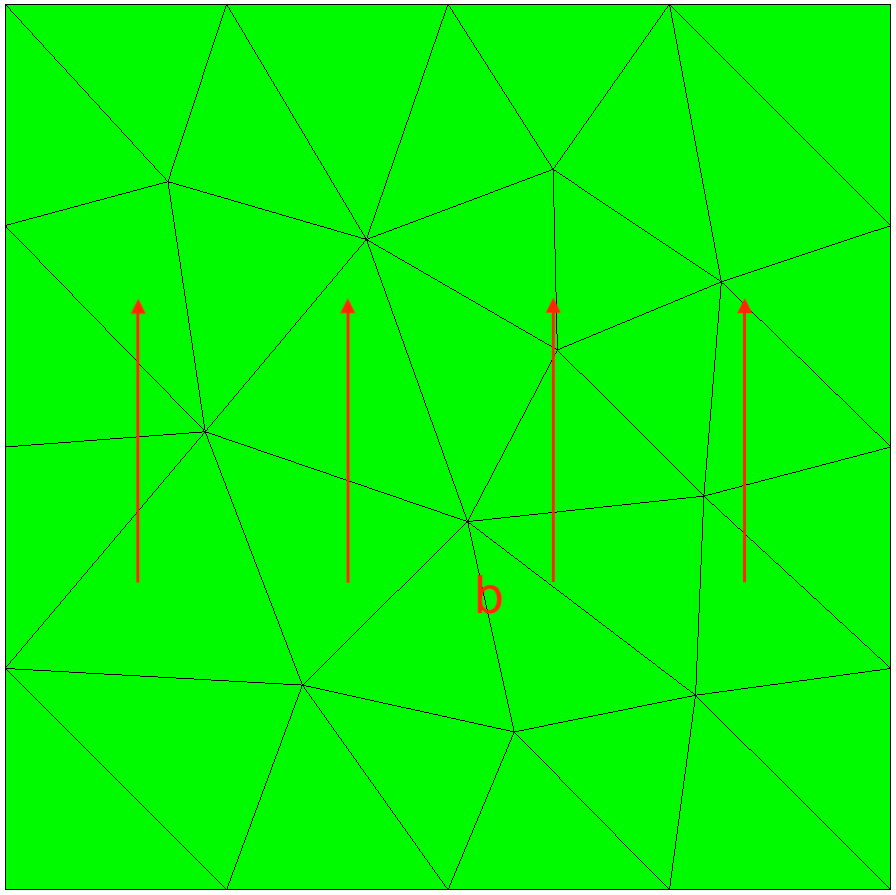}\\
(a) & (b)
\end{tabular}
\caption{Test~\ref{Num-Test2-1}: Example meshes for $h = 1/4$: (a) aligned quadrilateral mesh; (b) non-aligned triangular mesh.}\label{fig:TestSquare-1}
\end{figure}

\noindent\textbf{Problem Setting}: The computational domain is chosen as $\Omega = [0,1]\times[0,1]$, with a Dirichlet boundary condition at the vertical boundary $x = 1$, where the solution $u|_{x=1}$ is imposed. The rest of the boundaries are with homogeneous Neumann conditions. The magnetic field is vertical with components $b_1 = 0$ and $b_2 = 1$. The perpendicular diffusion is set as $D_{\perp} = 1$, and a source is imposed with the form
\begin{eqnarray*}
f = \sin(\pi x).
\end{eqnarray*}
The analytical solution is $u = \dfrac{\sin(\pi x)}{\pi^2}$. Since the value of the parallel diffusion does not affect the analytical solution, the effect on the numerical solution is entirely due to the numerical diffusion introduced by the scheme on a given discretization. Results from the aligned quadrilateral mesh (Figure~\ref{fig:TestSquare-1}a) and non-aligned triangular mesh (Figure~\ref{fig:TestSquare-1}b) are compared below.

\begin{figure}[H]
\centering
\begin{tabular}{cc}
\includegraphics[width =.45\textwidth]{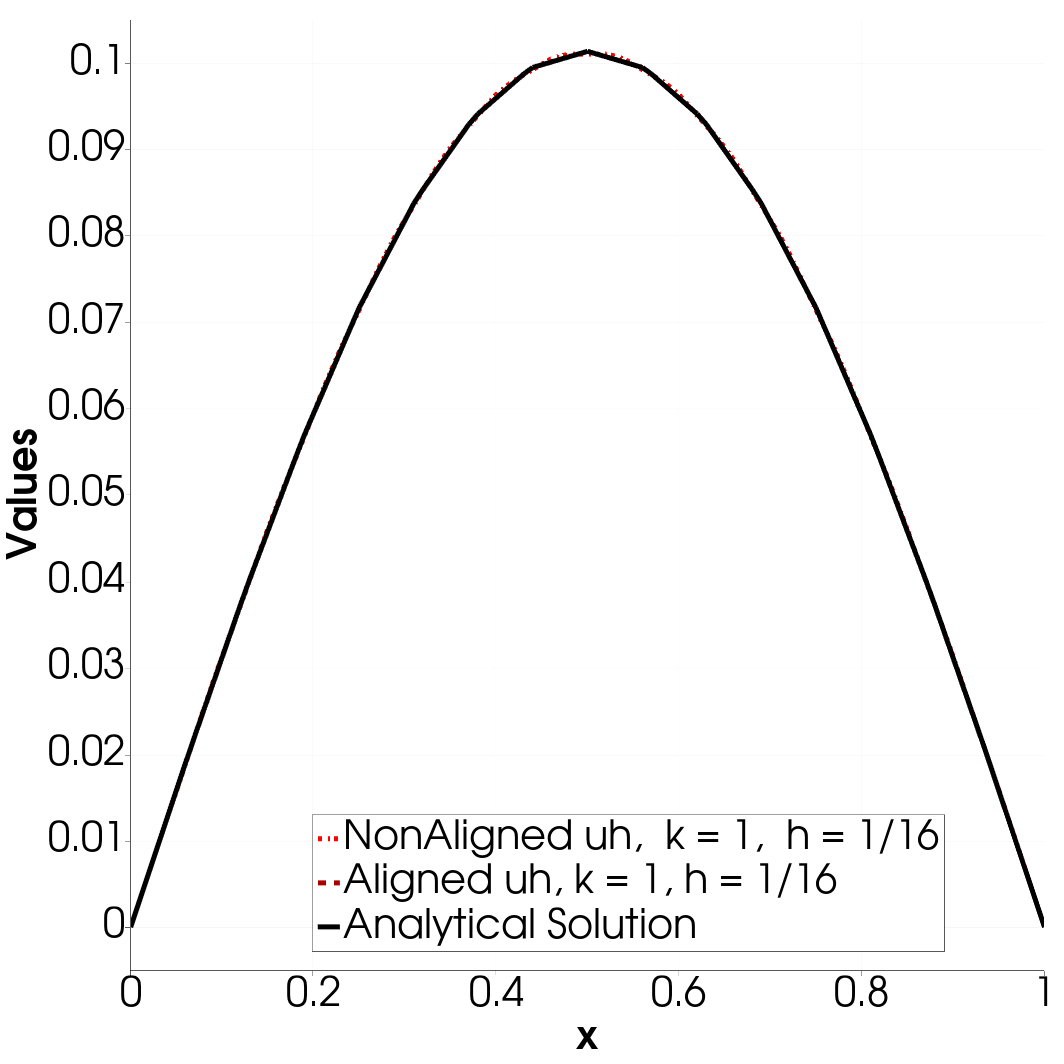}
&\includegraphics[width =.45\textwidth]{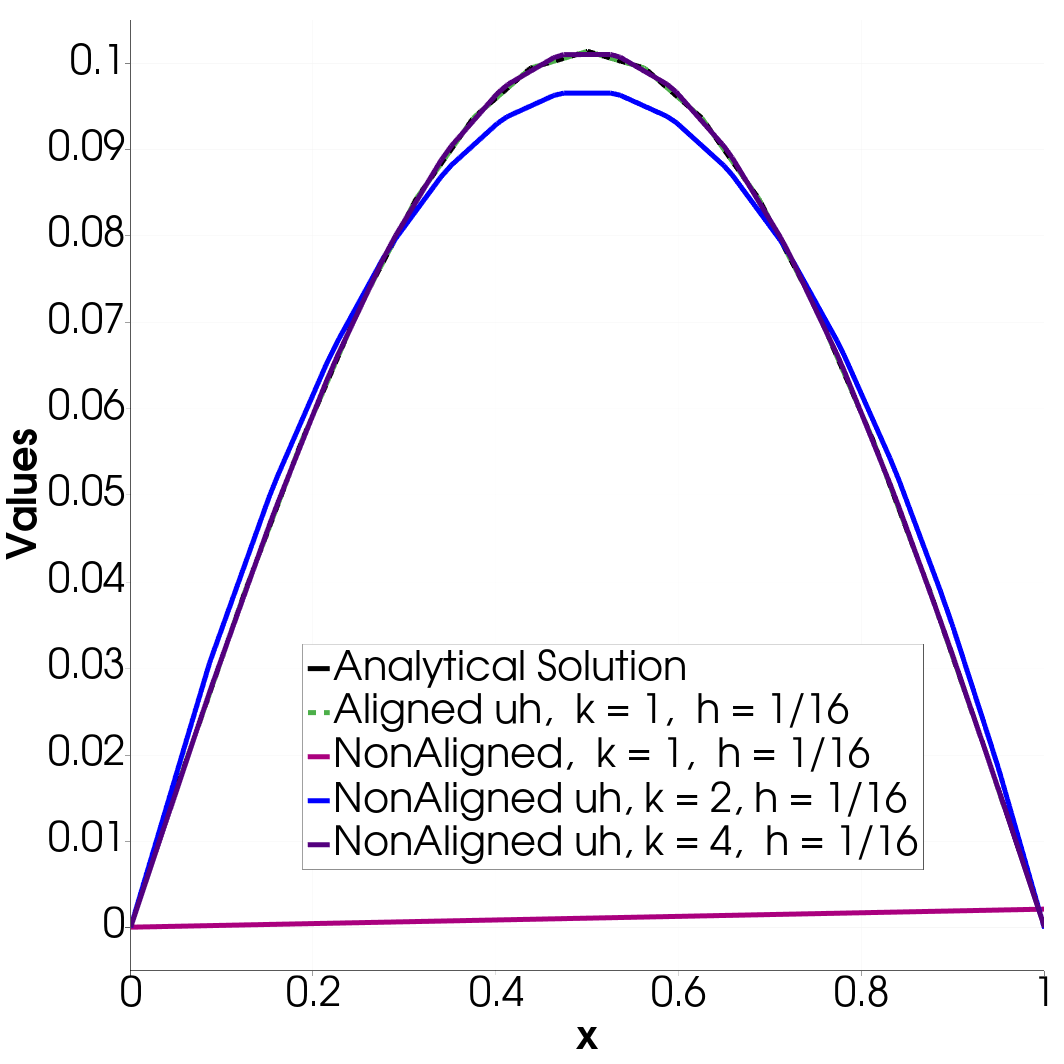}\\
(a) & (b)
\end{tabular}
\caption{Test~\ref{Num-Test2-1}: Results comparing aligned and non-aligned meshes with $h = 1/16$: (a) the isotropic case with $D_{\|} = 1$; (b) the anisotropic case with $D_{\|} = 10^9$.}\label{fig:TestSquare-2}
\end{figure}
We plot the numerical solution on the trace $y = 0.0625$ and compare it with the analytical solution in Figure~\ref{fig:TestSquare-2}. For the isotropic case with $D_{\|}= 1$, the linear ($k=1$) solutions on the non-aligned and aligned meshes agree with the analytical solution. In the strongly anisotropic case with $D_{\|} = 10^9$, the linear ($k=1$) numerical solution on the aligned mesh matches the exact solution. However, the low order numerical schemes on the non-aligned triangular mesh ($k = 1,2$) produce polluted numerical solutions. 
As we increase the polynomial degree up to $k = 4$, the numerical solution on the non-aligned triangular mesh matches the analytical solution. Thus, when the aligned mesh challenging to generate or use in a simulation, we may need to adopt high-order computational scheme to obtain satisfactory numerical solutions.

\subsubsection{Test: Diffusion of a Gaussian Source}\label{Num-Test2-2}
\begin{figure}[H]
\centering
\begin{tabular}{cc}
\includegraphics[width =.3\textwidth]{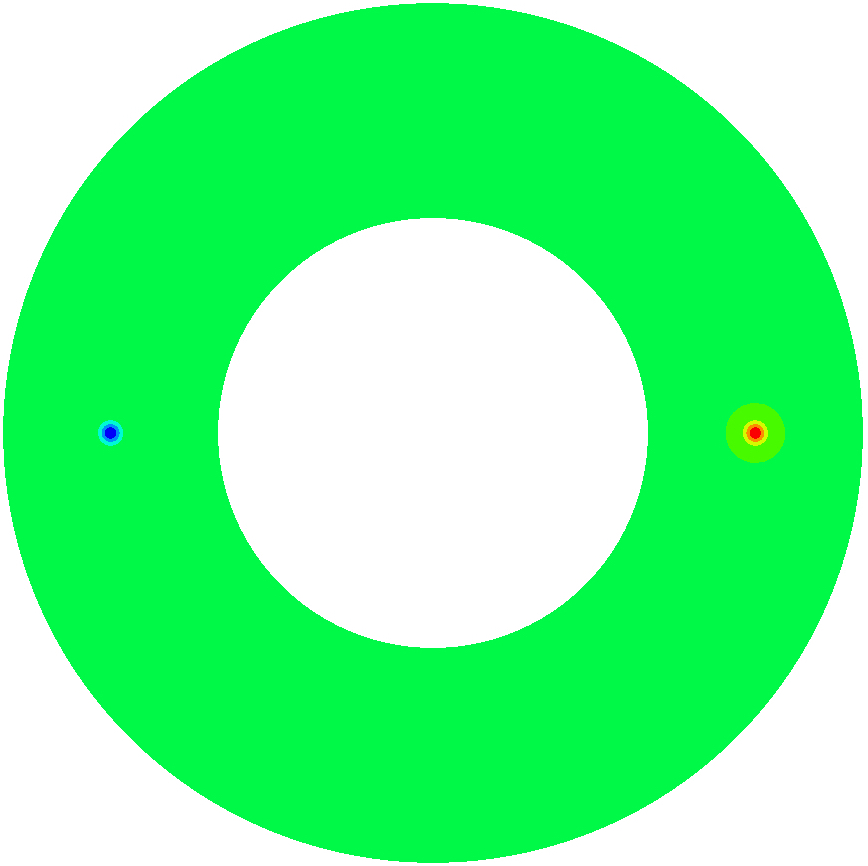}\quad\quad\quad\quad
&\includegraphics[width =.3\textwidth]{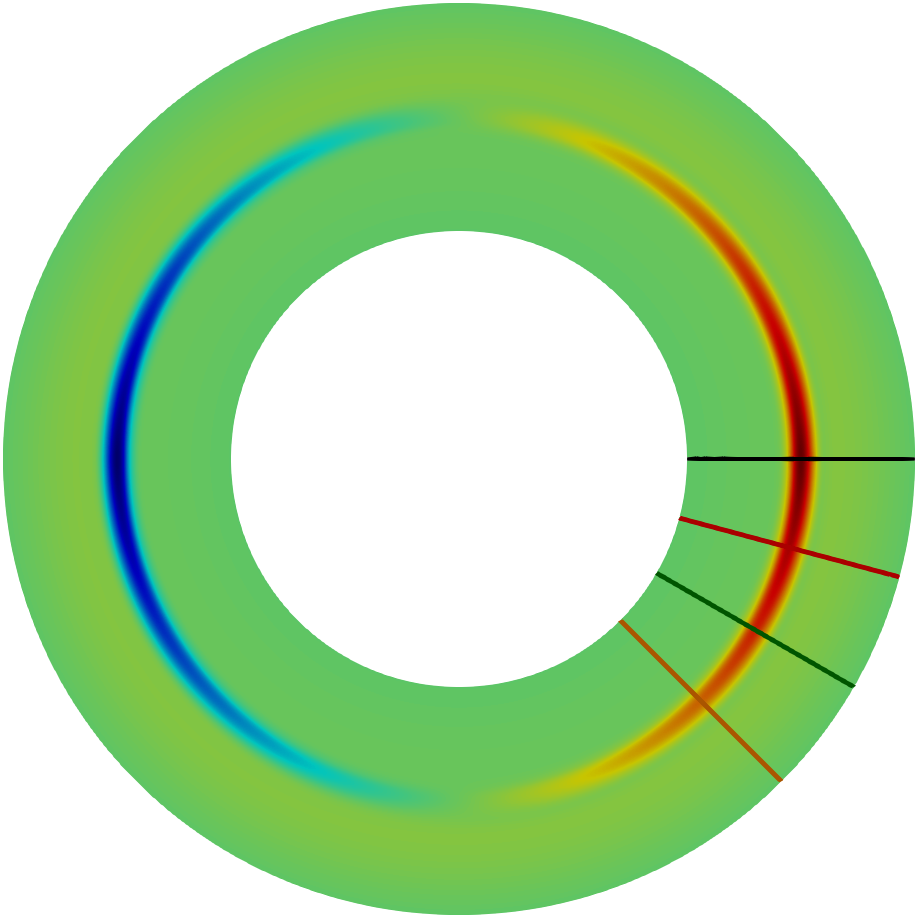} \\ 
(a) & (b) 
\end{tabular}
\caption{Test~\ref{Num-Test2-2}: Diffusion of a Gaussian Source with $D_{\|} = 10^9$: (a) 2D plot of the source; (b) Plot of the solution for $k = 8$ and $N_r = 8$.}\label{fig:GaussianSource}
\end{figure}
\noindent\textbf{Problem Setting}: 
Let the computational domain $\Omega$ be an annulus with exterior radius $r = 2$ and interior radius $r = 1$.
In this test, a Gaussian source is diffused towards an identical sink. The source and sink are defined as 
\begin{eqnarray*}
f_{sc} = D_{\|}\exp(-r_{sc}^2/0.05^2)\text{ and }f_{sk} = -D_{\|}\exp(-r_{sk}^2/0.05^2),
\end{eqnarray*}
where $r_{sc} = \sqrt{(x-1.5)^2+y^2}$ and $r_{sk} = \sqrt{(x+1.5)^2+y^2}$. The magnetic field direction is chosen as $b_1 = \dfrac{y}{r}, \ b_2 = -\dfrac{x}{r}$. In this test, we shall choose $D_{\|} = 10^9$ and $D_{\perp} = 1$. 

The source is plotted in Figure~\ref{fig:GaussianSource}a. The numerical solution for $k = 8$, $N_r = 8,N_\theta = 32$ on the quadrilateral aligned mesh is shown in Figure~\ref{fig:GaussianSource}b.

\begin{figure}[H]
\centering
\begin{tabular}{ccc}
\includegraphics[width =.3\textwidth]{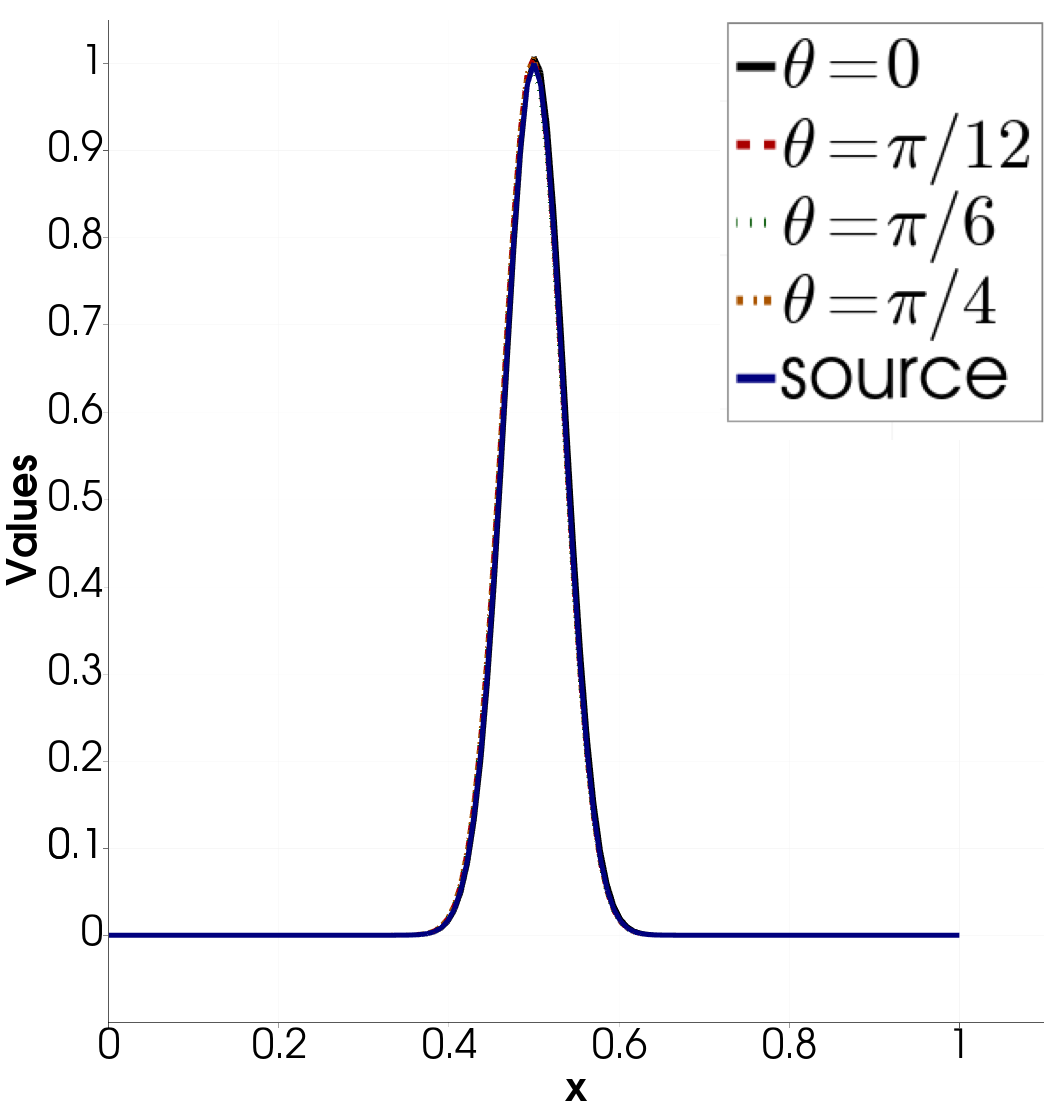}
&\includegraphics[width =.3\textwidth]{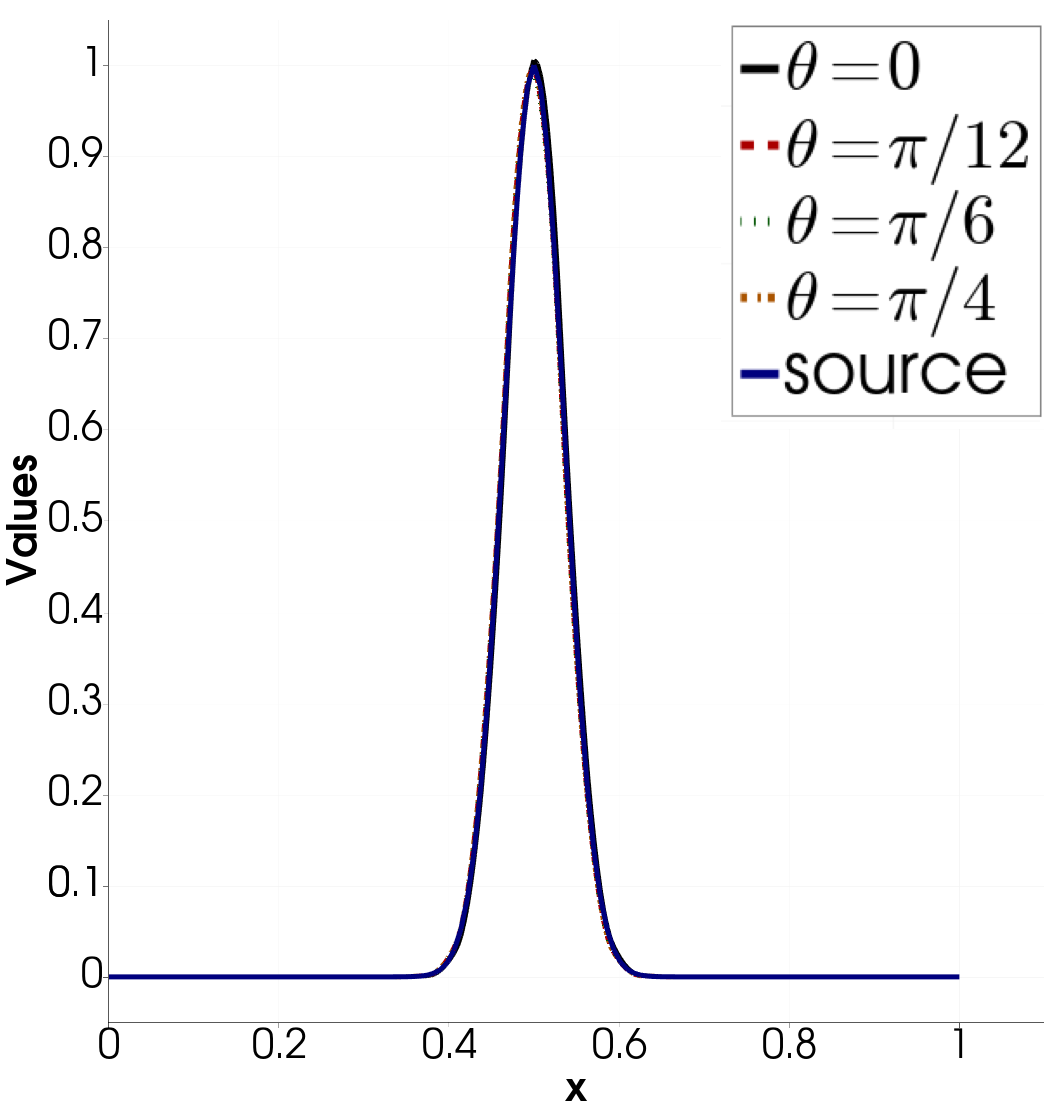}
& \includegraphics[width =.3\textwidth]{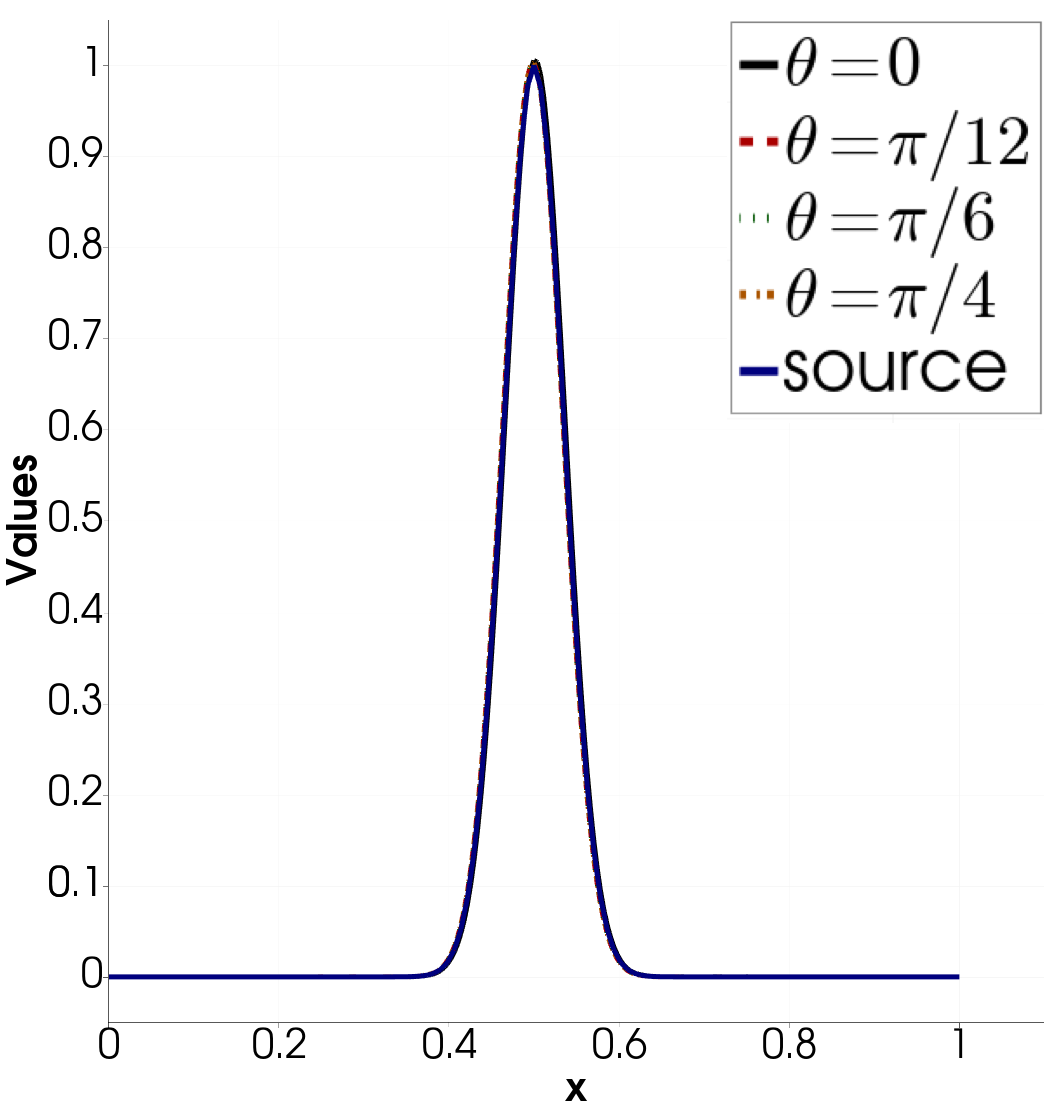}\\
(a) & (b) & (c)
\end{tabular}
\caption{Test~\ref{Num-Test2-2}: Trace plot with $D_{\|} = 10^9$ on the quadrilateral aligned mesh: (a) $k = 1$ and $N_r = 1/64$; (b) $k = 2$ and $N_r = 1/32$; (c) $k = 8$ and $N_r = 1/8$. The $x$-axis denotes
the span from the inner to the outer radius of the annulus.}
\label{fig:GaussianSource-1}
\end{figure}

Since $D_{\|}\gg D_{\perp}$, the solution is expected to be diffused only in the parallel direction. This means that the normalized profiles (with respect to the maximum on each radial line) should overlap the plot of $\exp(-((x-1.5)^2+y^2)/0.05^2)$. Therefore, the difference between normalized profiles and the expected profile is from the error introduced via polluted numerical diffusion. In Figure:~\ref{fig:GaussianSource-1} we show the trace plot along the four lines illustrated in Figure~\ref{fig:GaussianSource}b calculated using the quadrilateral aligned mesh for the angles $\theta = 0,\pi/12, \pi/6,\pi/4.$. For all the trace locations, the figures show that the aligned discretizations provide satisfied solutions at any polynomial degree $k$.

\begin{figure}[H]
\centering
\begin{tabular}{ccc}
\includegraphics[width =.3\textwidth]{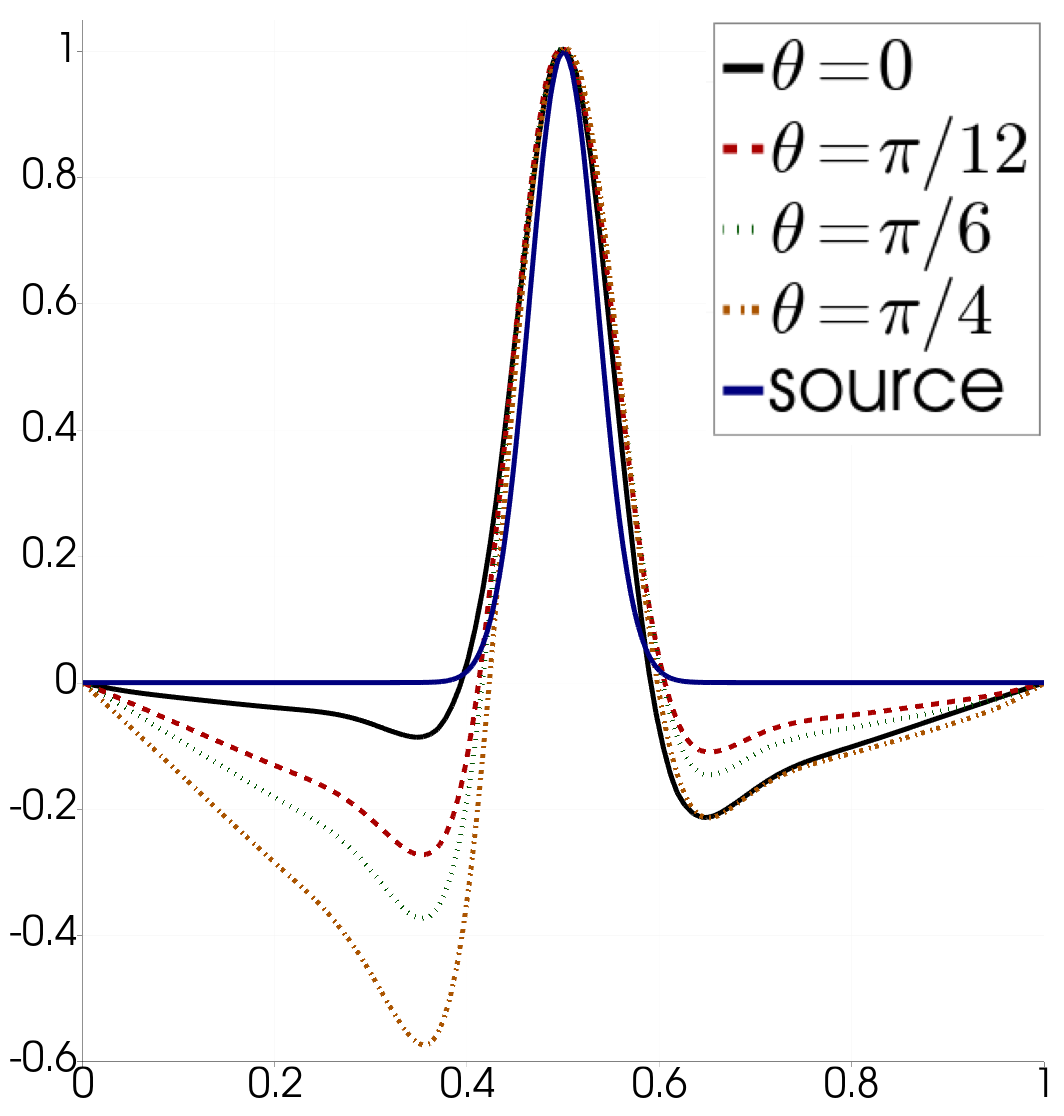}
&\includegraphics[width =.3\textwidth]{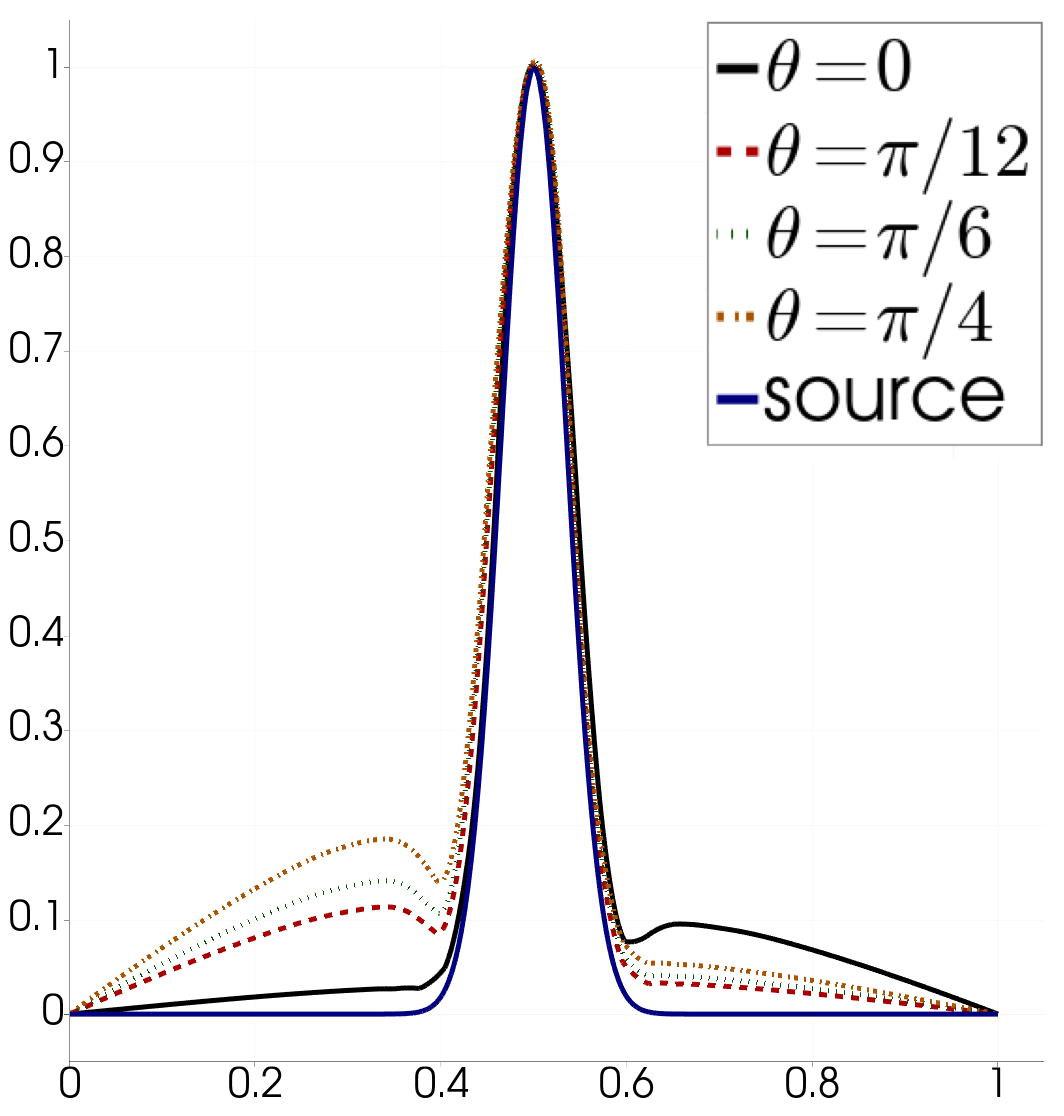}
& \includegraphics[width =.3\textwidth]{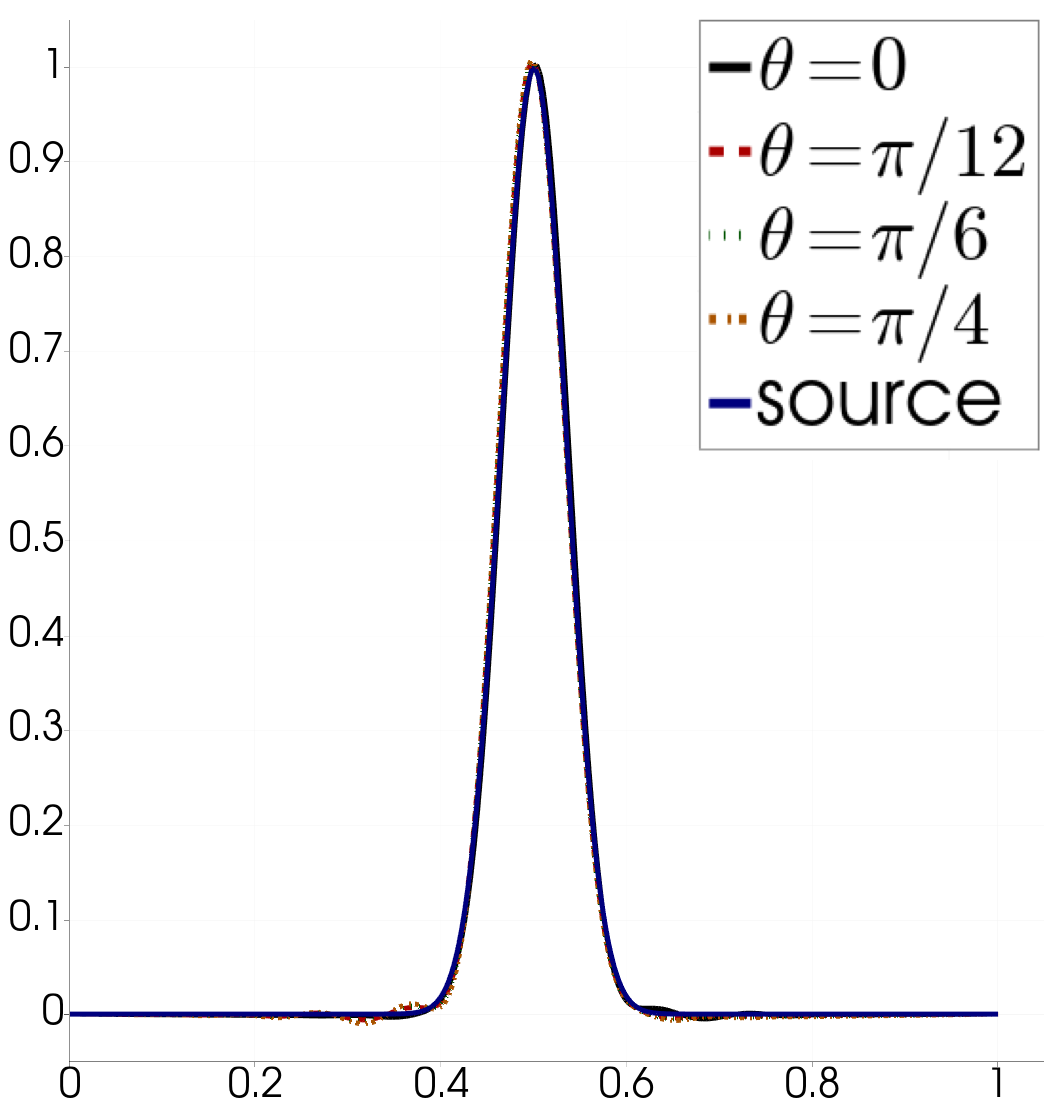}\\
(a) & (b) & (c)
\end{tabular}
\caption{Test~\ref{Num-Test2-2}: Trace plot with $D_{\|} = 10^9$ on the triangular non-aligned mesh: (a) $k = 1$ and $h = 1/64$; (b) $k = 2$ and $h = 1/32$; (c) $k = 8$ and $h = 1/8$. The $x$-axis denotes the span from the inner to the outer radius of the annulus. }\label{fig:GaussianSource-2}
\end{figure}

However, Figure~\ref{fig:GaussianSource-2} shows that the non-aligned triangular mesh does not behave well when the scheme order is low. As we can see from this figure, a significant spreading of the solution is visible for the low order elements ($k = 1$ and $k = 2$) but for high-order elements a satisfactory numerical solution is obtained. Again, it suggests we adopt a high-order scheme for the non-aligned mesh.

\subsubsection{Test: Two Magnetic Islands}\label{Num:Test2-3}
\begin{figure}[H]
\centering
\begin{tabular}{cc}
\includegraphics[width =.4\textwidth]{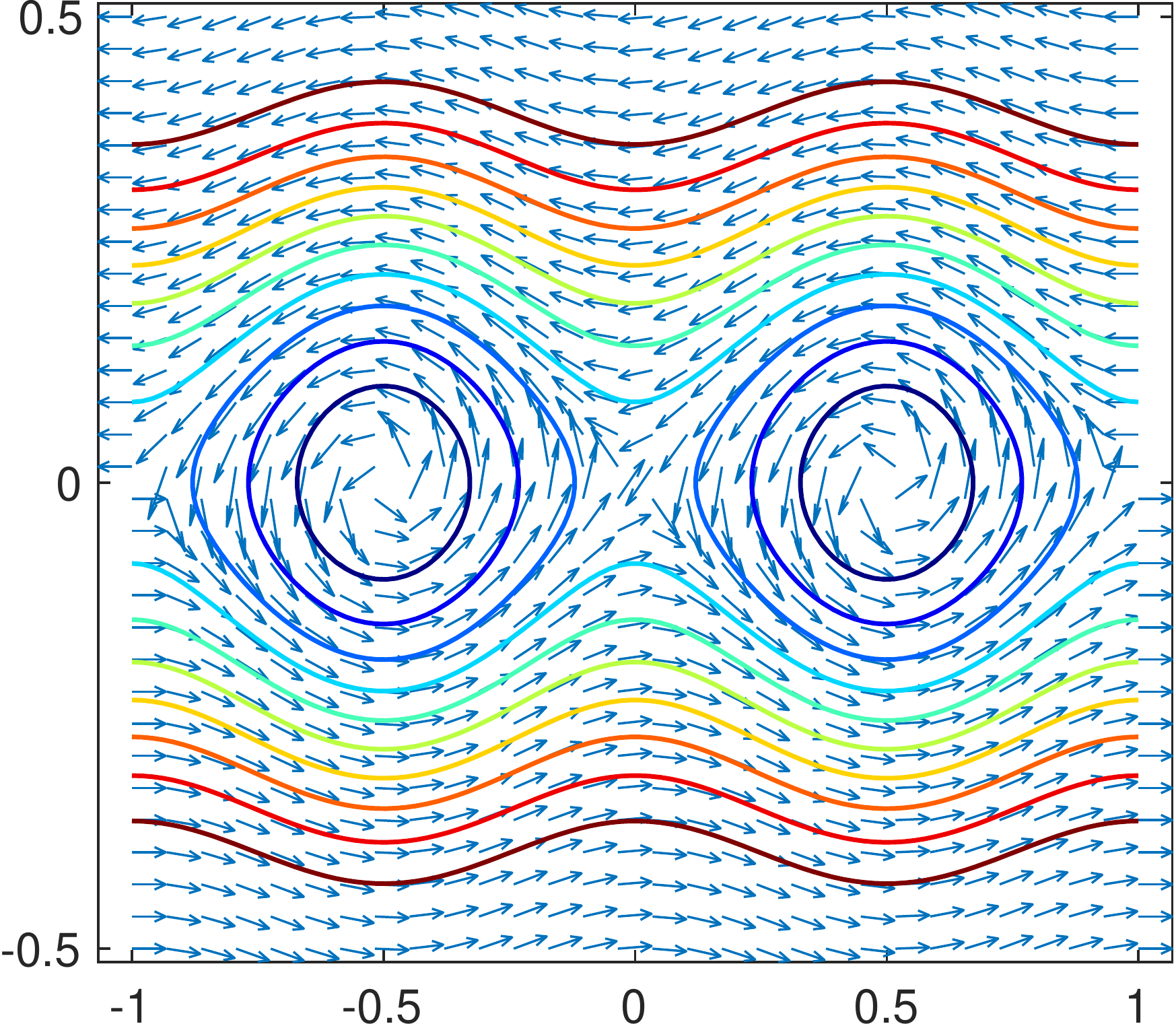}
&\includegraphics[width =.45\textwidth,height = .375\textwidth]{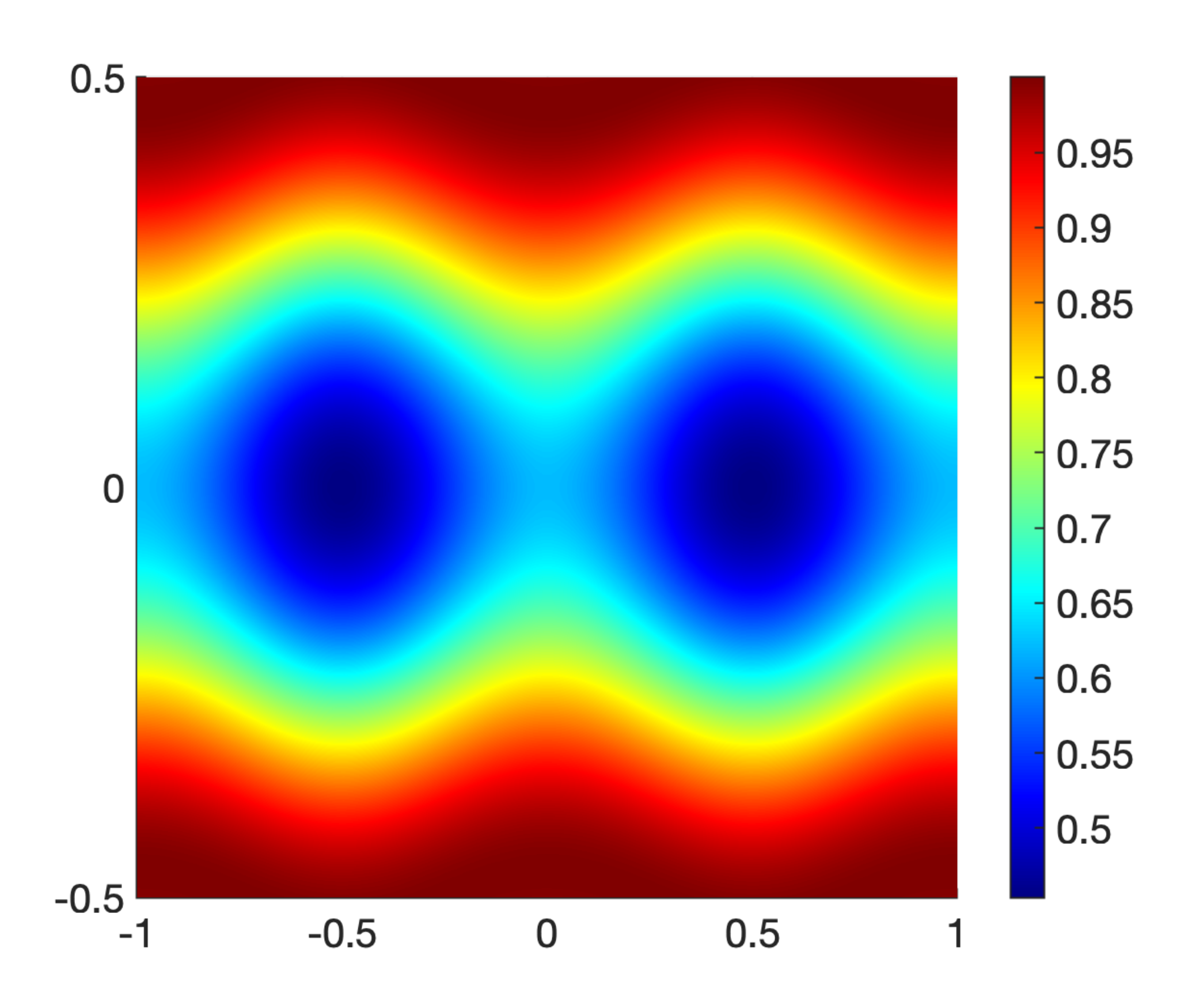}\\
(a) & (b)
\end{tabular}
\caption{Test~\ref{Num:Test2-3}: Plot of (a) magnetic field ${\bf b}$; (b) exact solution.}\label{fig:NumTest-2-3-1}
\end{figure}
\noindent\textbf{Problem Setting}: 
Let the computational domain be $\Omega = [-1,1]\times[-0.5,0.5]$ and the exact solution be
\begin{eqnarray}
u= \cos(\frac{1}{10}\cos(2\pi(x-3/2))+\cos(\pi y)). 
\end{eqnarray}
Let ${\bf b} = (b_1,b_2)^\top$ and
\begin{eqnarray}
{\bf b} = \frac{\bf B}{\bf |B|},\ {\bf B} = \begin{pmatrix}
-\pi\sin(\pi y)\\ \frac{2\pi}{10} \sin(2\pi(x-3/2))
\end{pmatrix}.
\end{eqnarray}
The perpendicular diffusion is set as $D_{\perp}= 1$ with varying values of $D_\|$ for the diffusion tensor (\ref{eq:diff-coef}).
The magnetic field and the exact solution are plotted in Figure~\ref{fig:NumTest-2-3-1}. It is noted that the magnetic field shows two islands. 

Due to the complexity of the magnetic direction, it is almost impossible to generate an aligned mesh for resolving the problem's anisotropy. Here we shall investigate the performance of the high order scheme on the Cartesian grid with $1/h$ partitions in both $x$ and $y$ directions.

\begin{table}[H]
 \centering
     \tabcolsep=3pt
     {\footnotesize
\begin{tabular}{c||cc|cc|cc|cc|cc}\hline\hline
& \multicolumn{2}{c|}{$D_{\|}$ = 1E1} & \multicolumn{2}{c|}{$D_{\|}$ = 1E2}  & \multicolumn{2}{c|}{$D_{\|}$ = 1E4}  & \multicolumn{2}{c|}{$D_{\|}$ = 1E6} & \multicolumn{2}{c}{$D_{\|}$ = 1E8} \\ 
$1/h$ & $\|u-u_h\|$ & order & $\|u-u_h\|$ & order & $\|u-u_h\|$ & order & $\|u-u_h\|$ & order & $\|u-u_h\|$ & order\\ \hline\hline
\multicolumn{11}{c}{$k = 1$}\\
8	&2.44E-02	&		&5.54E-02	&		&7.12E-02	&		&7.14E-02	&		&7.14E-02 & \\	
16	&8.31E-03	&1.56	&3.26E-02	&0.77	&5.66E-02	&0.33	&5.70E-02	&0.33	&5.70E-02	&0.33\\
32	&2.39E-03	&1.80	&1.40E-02	&1.22	&4.73E-02	&0.26	&4.85E-02	&0.23	&4.85E-02	&0.23\\
64	&6.33E-04	&1.92	&4.50E-03	&1.64	&3.90E-02	&0.28	&4.27E-02	&0.18	&4.28E-02	&0.18\\
128	&1.62E-04	&1.97	&1.24E-03	&1.86	&2.84E-02	&0.46	&3.86E-02	&0.15	&3.87E-02	&0.14\\
256	&4.10E-05	&1.99	&3.23E-04	&1.94	&1.51E-02	&0.91	&3.54E-02	&0.12	&3.59E-02	&0.11\\ \hline\hline
\multicolumn{11}{c}{$k = 2$}\\						
4	&5.68E-03	&		&1.20E-02	&		&1.73E-02	&		&1.74E-02	&		&1.74E-02 &\\	
8	&6.01E-04	&3.24	&1.36E-03	&3.14	&7.80E-03	&1.15	&8.26E-03	&1.08	&8.26E-03	&1.07\\
16	&7.16E-05	&3.07	&1.22E-04	&3.49	&3.22E-03	&1.28	&4.91E-03	&0.75	&4.94E-03	&0.74\\
32	&8.84E-06	&3.02	&1.11E-05		&3.45	&4.93E-04	&2.71	&3.35E-03	&0.55	&3.56E-03	&0.47\\
64	&1.10E-06	&3.00	&1.19E-06	&3.23	&3.67E-05	&3.75	&1.56E-03	&1.11	&2.85E-03	&0.32\\
128	&1.38E-07	&3.00	&1.40E-07	&3.08	&2.42E-06	&3.92	&1.95E-04	&3.00	&2.23E-03	&0.35\\ \hline\hline
\multicolumn{11}{c}{$k = 3$}\\				
4	&7.04E-04	&		&1.15E-03	&		&6.65E-03	&		&7.34E-03	&		&7.35E-03 &\\	
8	&4.48E-05	&3.98	&6.43E-05	&4.16	&5.26E-04	&3.66	&9.12E-04	&3.01	&9.41E-04	&2.96\\
16	&2.78E-06	&4.01	&3.16E-06	&4.35	&5.17E-05	&3.35	&1.05E-03	&-0.21	&1.33E-03	&-0.50\\
32	&1.74E-07	&4.00	&1.80E-07	&4.13	&1.36E-06	&5.25	&8.98E-05	&3.55	&8.84E-04	&0.59\\
64	&1.09E-08	&4.00	&1.10E-08	&4.03	&3.79E-08	&5.16	&2.32E-06	&5.27	&1.32E-04	&2.74\\
128	&6.83E-10	&4.00	&6.83E-10	&4.01	&1.25E-09	&4.92	&6.45E-08	&5.17	&3.37E-06	&5.29\\ \hline\hline
\multicolumn{11}{c}{$k = 4$}\\			
4	&3.64E-05	&		&5.19E-05	&		&9.96E-04	&		&3.08E-03	&		&3.14E-03&\\	
8	&1.03E-06	&5.15	&1.50E-06	&5.11	&1.70E-05	&5.88	&9.21E-05	&5.06	&9.92E-05	&4.99\\
16	&3.04E-08	&5.08	&3.76E-08	&5.32	&1.34E-07	&6.99	&1.31E-06	&6.14	&8.19E-06	&3.60\\
32	&9.16E-10	&5.05	&1.02E-09	&5.20	&2.45E-09	&5.77	&9.02E-09	&7.18	&5.38E-07	&3.93\\
64	&2.81E-11		&5.03	&2.96E-11		&5.11	&5.56E-11		&5.46	&4.43E-10	&4.35	&6.52E-08	&3.04\\ \hline\hline
\multicolumn{11}{c}{$k = 5$}\\
2	&2.62E-04	&		&3.71E-04	&		&4.12E-03	&		&5.22E-03	&		&5.23E-03&\\	
4	&3.12E-06	&6.39	&4.81E-06	&6.27	&9.10E-05	&5.50	&3.43E-04	&3.93	&3.53E-04	&3.89\\
8	&7.58E-08	&5.36	&8.81E-08	&5.77	&5.14E-07	&7.47	&1.59E-05	&4.43	&9.97E-05	&1.82\\
16	&1.09E-09	&6.12	&1.14E-09	&6.27	&5.62E-09	&6.52	&9.92E-08	&7.32	&2.01E-06	&5.63\\
32	&1.66E-11		&6.03	&1.68E-11		&6.09	&3.92E-11		&7.16	&4.08E-10	&7.93	&1.42E-08	&7.15\\ \hline\hline
\multicolumn{11}{c}{$k = 6$}\\										
1	&1.50E-03	&		&3.18E-03	&		&2.83E-02	&		&5.24E-02	&		&5.29E-02&\\	
2	&1.70E-05	&6.46	&3.72E-05	&6.42	&7.55E-04	&5.23	&1.35E-03	&5.28	&1.36E-03	&5.28\\
4	&3.18E-07	&5.74	&5.40E-07	&6.11	&5.49E-06	&7.10	&1.10E-04	&3.61	&1.73E-04	&2.97\\
8	&1.91E-09	&7.38	&2.67E-09	&7.66	&2.84E-08	&7.60	&5.60E-07	&7.62	&9.63E-06	&4.17\\
16	&1.37E-11		&7.12	&1.63E-11		&7.36	&7.90E-11		&8.49	&2.53E-09	&7.79	&1.57E-08	&9.26\\ \hline\hline
\end{tabular}
}
\caption{Test~\ref{Num:Test2-3}: Error profiles and convergence results.}\label{tab:NumTest-2-3-1}
\end{table}


Table~\ref{tab:NumTest-2-3-1} shows the error profiles and the convergence results for varying values of $D_{\|}$ and polynomial degrees $k$ by using the uniform rectangular meshes. It shows that with small anisotropy ($D_{\|} = $1 and 1E2), the scheme can still produce the optimal rate of convergence, which is at the order $\mathcal{O}(h^{k+1})$.  As expected, increasing the values of $D_{\|}$ also increases the $L^2$-error. However, for low order schemes ($k = $1 and 2) we do not observe the correct convergence rates. When the mesh is fine enough, the use of cubic elements is able to produce satisfactory numerical solutions with the numerical errors at the order below $10^{-6}$. In addition, one can observe that by increasing the polynomials' degree further, we obtain the desired numerical results with the extreme anisotropic ratios $D_{\|}/D_{\perp}$.

\subsection{Test: Convergence Test for Aligned and Non-aligned Meshes}\label{Num-Test3}
\begin{figure}[H]
\centering
\begin{tabular}{cc}
\includegraphics[width =.4\textwidth]{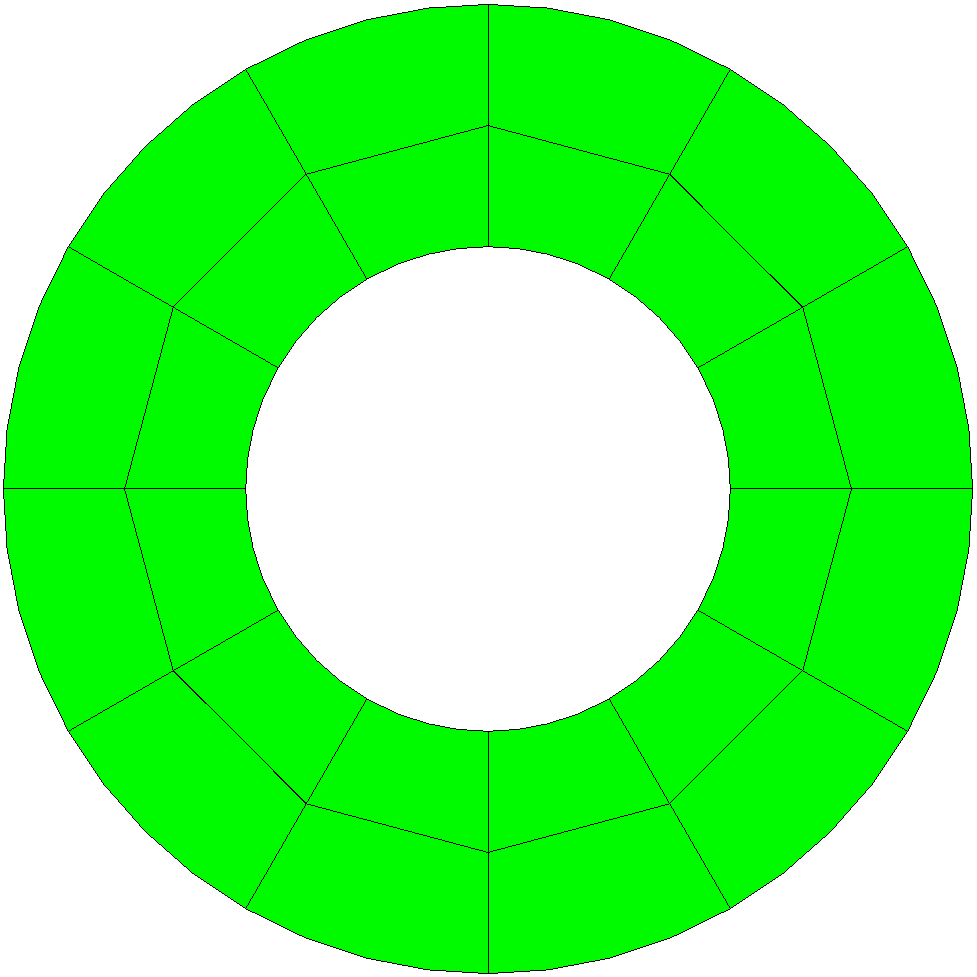}
&\includegraphics[width =.4\textwidth]{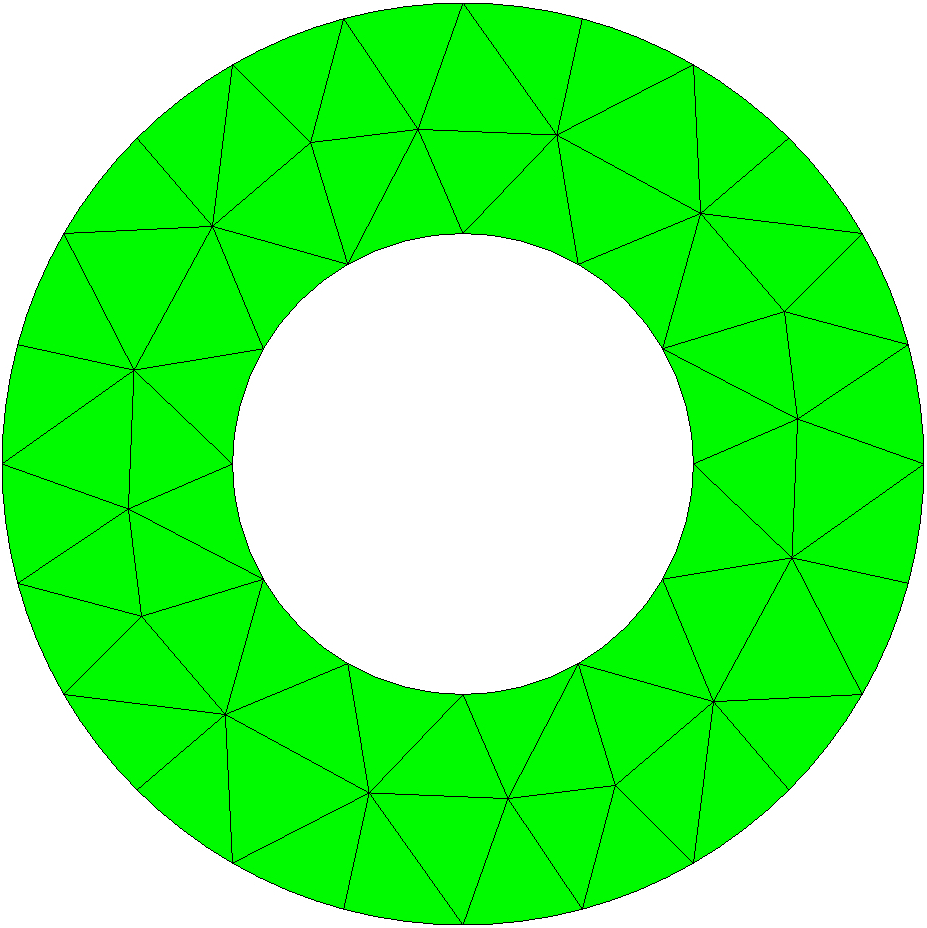}\\
(a) & (b)
\end{tabular}
\caption{Test~\ref{Num-Test3}: computational mesh with $N_r = 2$: (a) quadrilateral aligned mesh; (b) triangular non-aligned mesh.}\label{fig:TestSquare-mesh}
\end{figure}

In this test, we shall choose the manufactured solutions in (\ref{ex:TestAnnulus}) to test the numerical performance on the aligned quadrilateral mesh and non-aligned triangular meshes (as illustrated in Figure~\ref{fig:TestSquare-mesh}). 

\begin{figure}[H]
\centering
\begin{tabular}{cc}
\includegraphics[width =.48\textwidth]{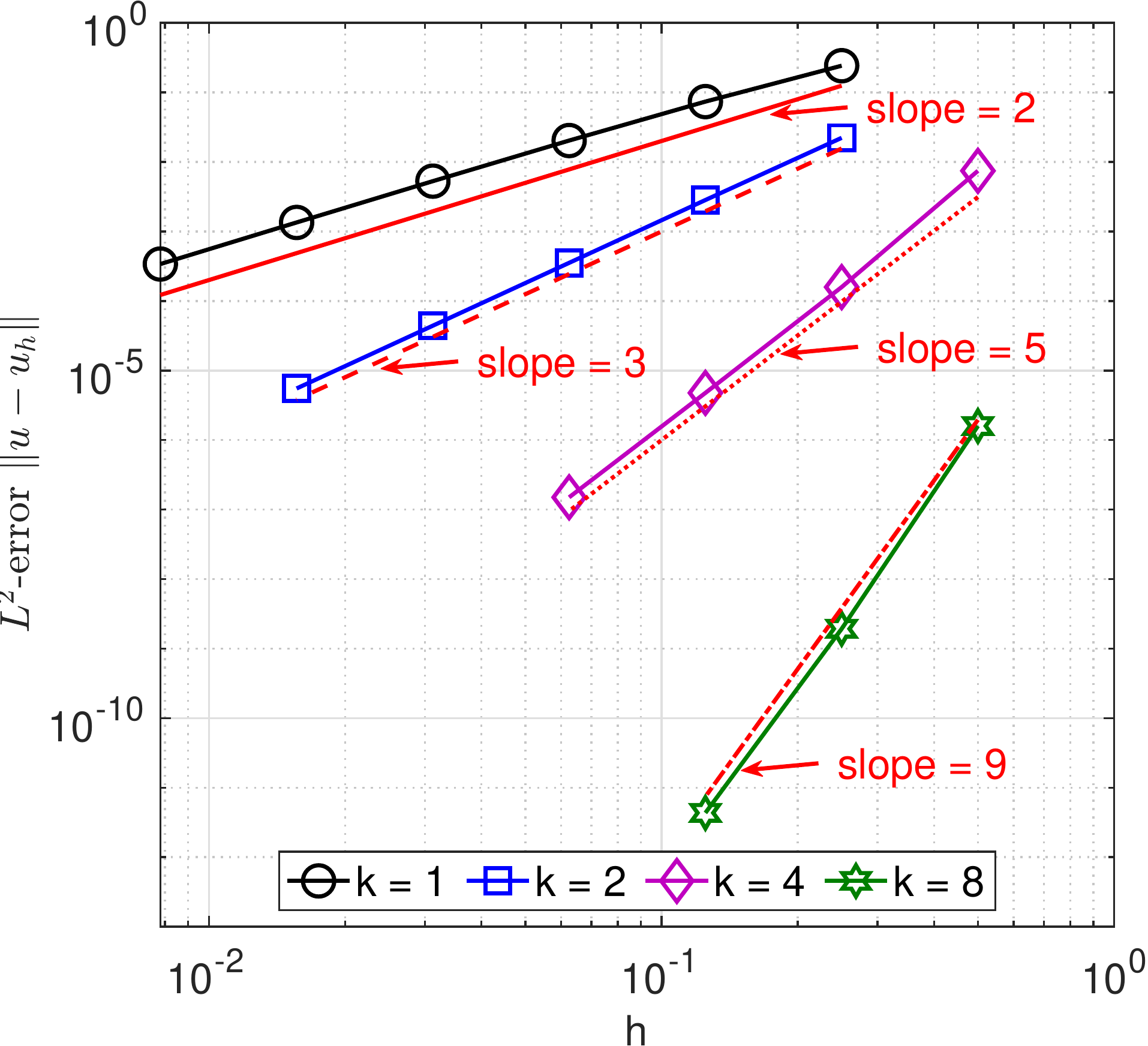}
&\includegraphics[width =.48\textwidth]{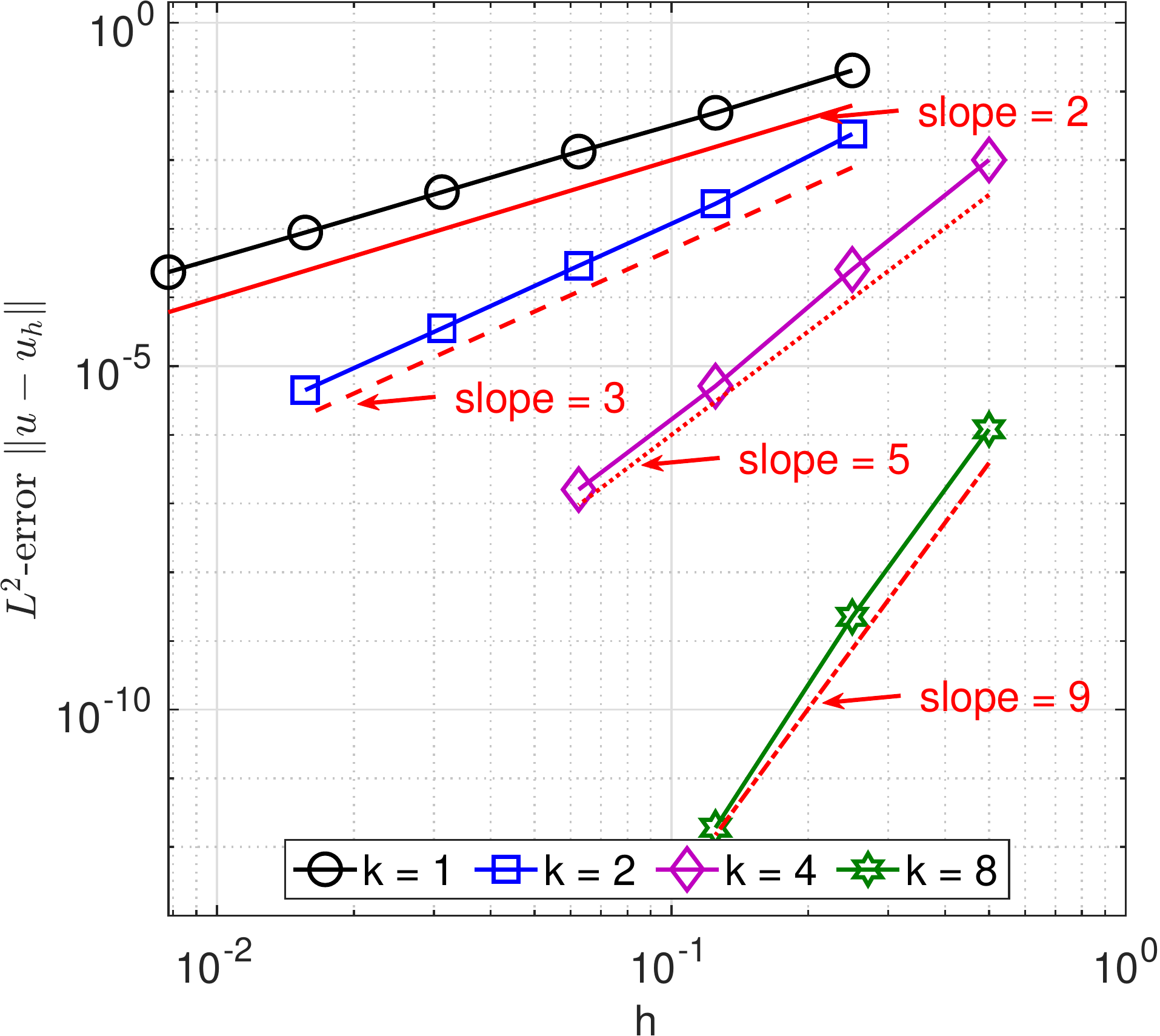}\\
(a) & (b)
\end{tabular}
\caption{Test~\ref{Num-Test3}: results of $L^2$-error convergence test on the meshes with $D_{\|} = 1$: (a) quadrilateral aligned mesh; (b) triangular non-aligned mesh.}\label{fig:TestSquare-3}
\end{figure}

\begin{figure}[H]
\centering
\begin{tabular}{cc}
\includegraphics[width =.48\textwidth]{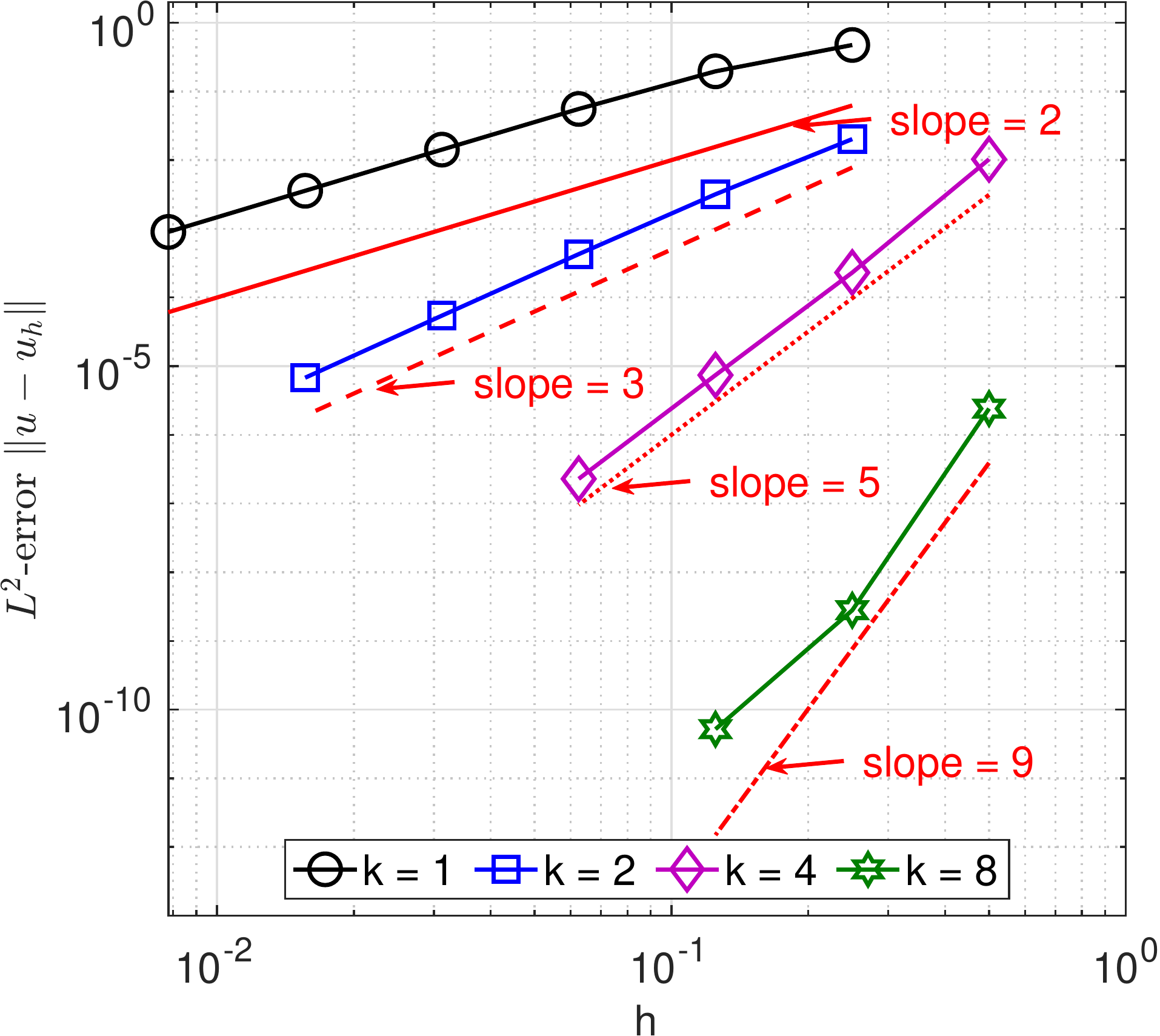}
&\includegraphics[width =.48\textwidth]{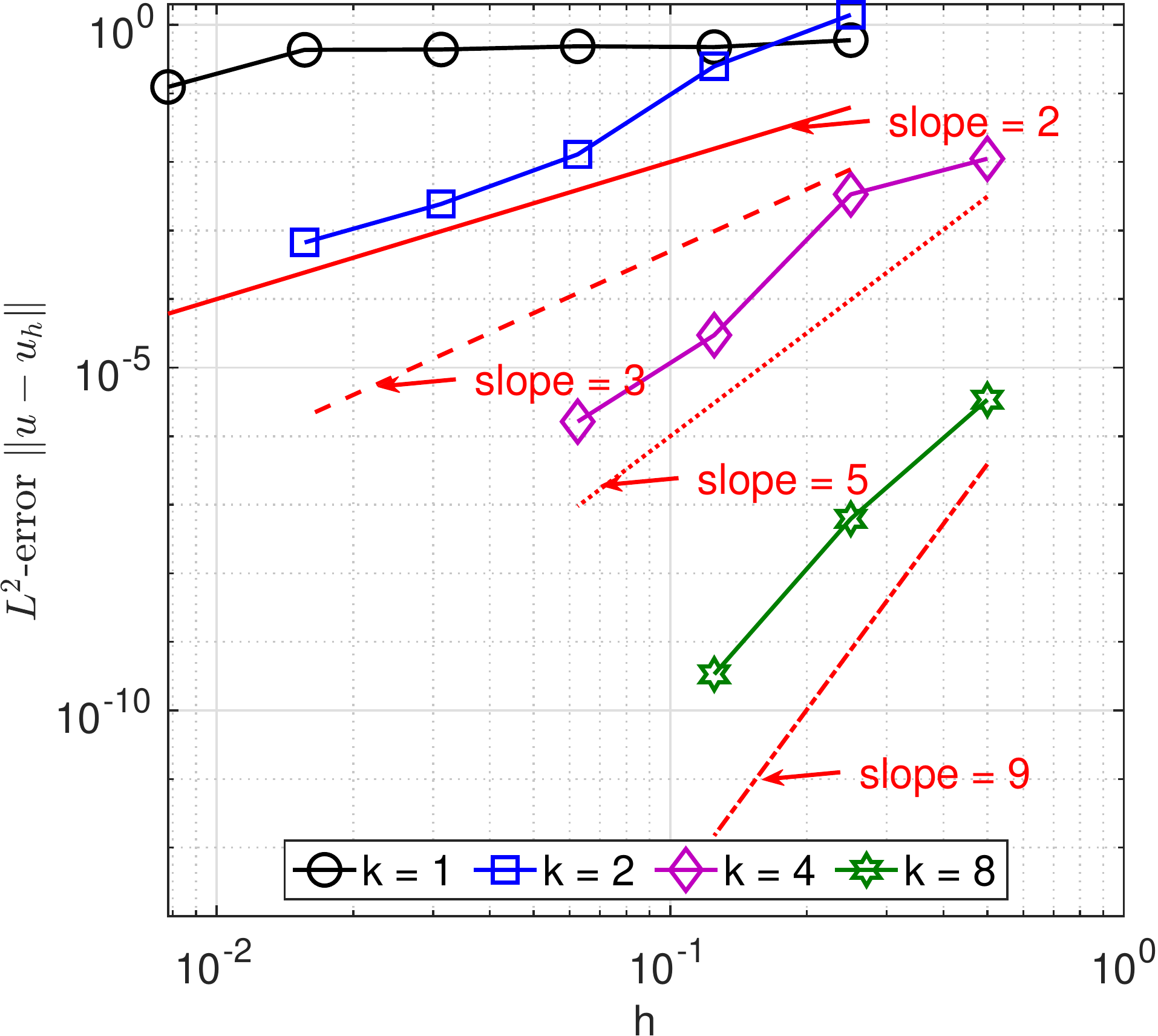}\\
(a) & (b)
\end{tabular}
\caption{Test~\ref{Num-Test3}: results of $L^2$-error convergence test on the meshes with $D_{\|} = 10^6$: (a) quadrilateral aligned mesh; (b) triangular non-aligned mesh.}\label{fig:TestSquare-4}
\end{figure}

Results are shown in Figure~\ref{fig:TestSquare-3} and Figure~\ref{fig:TestSquare-4} for an isotropic case ($D_{\|} = 1$) and an anisotropic case ($D_{\|} = 10^6$), respectively. We summarize our conclusions as below.

For the isotropic case (Figure~\ref{fig:TestSquare-3}), the convergence of $L^2$-error is nearly independent of the features of the mesh. The theoretical convergence rate $\mathcal{O}(h^{k+1})$ is achieved for all the polynomial approximations.

For the strongly anisotropic case, the convergence results of the $L^2$-error depends on the nature of the mesh type for low order elements.  In the aligned quadrilateral mesh, the optimal rate of convergence can still be obtained, as indicated by the slope plots in Figure~\ref{fig:TestSquare-4}a.  In the triangular mesh, by using linear elements, we barely observe any convergence order with respect to mesh size $h$.  In the case of quadratic elements, we can can obtain second order convergence but it should be third order. When the polynomial's degree is increased up to $k = 4$ and $k = 8$, the correct convergence results are recovered by the simulation. 
As such, the general conclusion of this section is that as we increase the polynomial degree, the impact on the non-alignment of the mesh is reduced and we recover the expected convergence behavior.

\subsection{Conditioning of the Corresponding Linear System}\label{Num-ConTest}
This sub-section is contributed to investigate the conditioning results of the derived IPDG schemes in the strong anisotropy cases. The problem setting is chosen the same as Test~\ref{Num-Test1}. In the following, we shall  demonstrate the performance (in terms of condition number) of the aligned and non-aligned meshes.

\begin{figure}[H]
\centering
\begin{tabular}{cc}
\includegraphics[width =.45\textwidth]{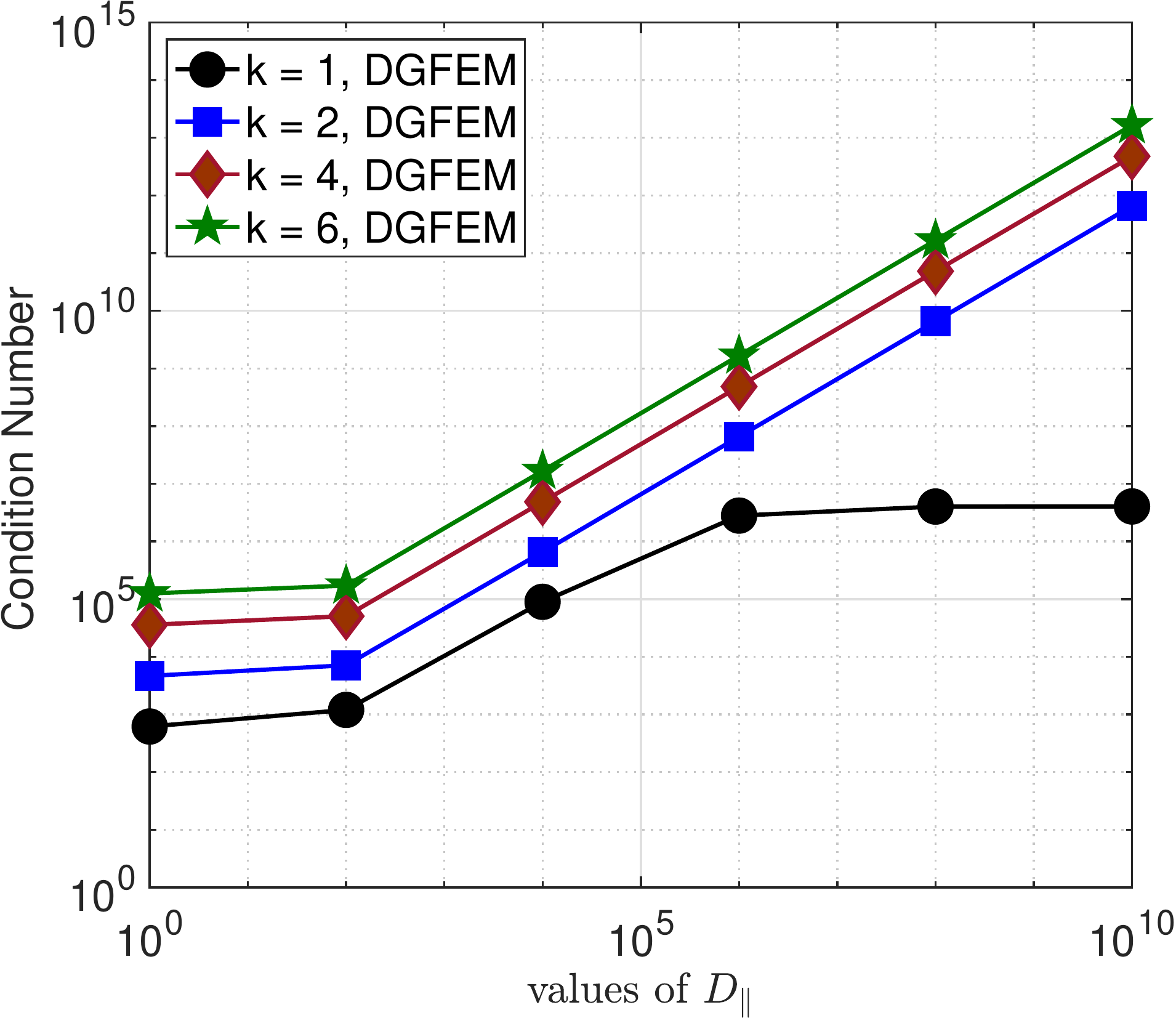}
&\includegraphics[width =.45\textwidth]{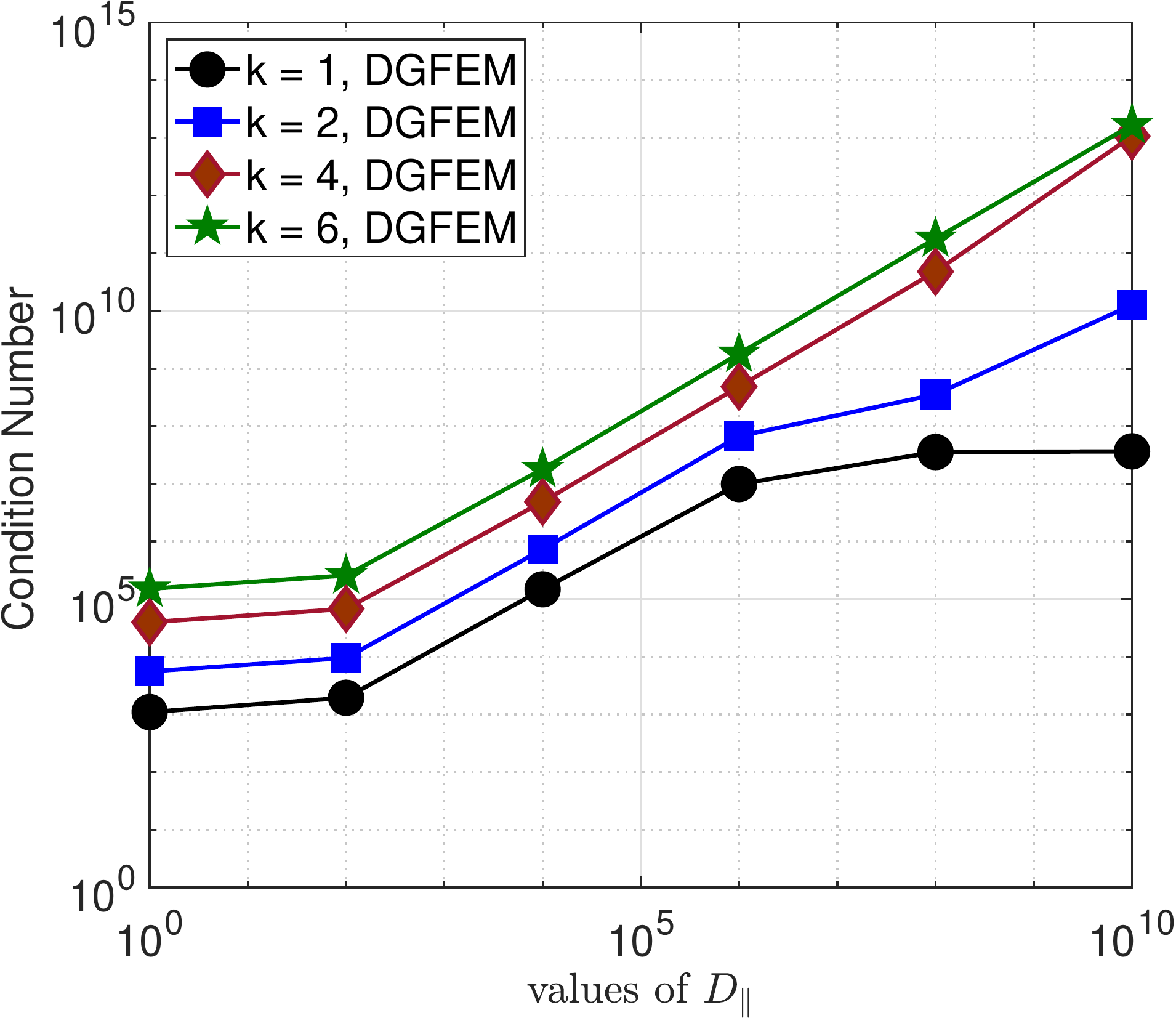}
\\
(a) & (b)
\end{tabular}
\caption{Test~\ref{Num-ConTest}. Plot of condition numbers with respect to  $D_{\|}$ on the aligned (a) Figure~\ref{fig:TestSquare-mesh}a; and non-aligned (b) Figure~\ref{fig:TestSquare-mesh}b meshes.}\label{fig:cond-L2_Quad}
\end{figure}

In Figure~\ref{fig:cond-L2_Quad} we test the performance on the quadrilateral and triangular meshes (from Figure~\ref{fig:TestSquare-mesh}) in the condition number sense for accuracy. The performance is very similar for both mesh types. The computed condition numbers (computed by the command ``condest" in PyAMG). It is noted that the conditional number is increasing almost at the order of $D_{\|}$ except when $k = 1$ where it appears bounded as we keep increasing the values in $D_{\|}$. We are still not clear about this phenomenon and will leave it for interested readers to investigate. Thus in the scheme with $k\ge 2$, the derived linear systems may have extremely large condition numbers in the strong anisotropic diffusion case. In such cases, an effective preconditioner is needed for an accurate and efficient numerical solver.





%
%
\section{Auxiliary Space Preconditioners}\label{sect:preconditioner}
In this section, we introduce auxiliary space preconditioners (ASP)~\cite{N1992,X1996} for solving the IPDG discretization~\eqref{eqn:DG}.  The IPDG method~\eqref{eqn:DG} can be written in the following operator form: given $\bm{f} \in V_{\text{DG}}'$, find $\bm{u} \in V_{\text{DG}}$ such that
\begin{equation*}
	\bm{A}_{\text{DG}}\bm{u} = \bm{f}, 
\end{equation*}
where $V_{\text{DG}}'$ is the dual of $V_{\text{DG}}$ and $\bm{A}: V_{\text{DG}} \mapsto V_{\text{DG}}'$ is the operator corresponding to the bilinear form $A_{\text{DG}}(\cdot, \cdot)$.  In order to design preconditioners using the auxiliary space preconditioning framework, as suggested in~\cite{ASVZ2017}, we consider the following product of auxiliary spaces:
\begin{equation*}
\bar{V}_{\text{DG}} = V_{\text{DG}} \times V_{\text{CG}}.
\end{equation*}
Note that, since $V_{\text{CG}} \subset V_{\text{DG}}$, this product of the auxiliary spaces can also be considered as a special subspace decomposition of $V_{\text{DG}}$, i.e., $V_{\text{DG}} = V_{\text{DG}} + V_{\text{CG}}$. Naturally, $V_{\text{CG}}$ is endowed with inner product $A_{\text{CG}}(\cdot, \cdot)$ whose operator form is $\bm{A}_{\text{CG}}: V_{\text{CG}} \mapsto V_{\text{CG}}'$ with $V_{\text{CG}}'$ being the dual of $V_{\text{CG}}$.  In addition, the transfor operator $\bm{\Pi}: V_{\text{CG}} \mapsto V_{\text{DG}}$ is the natrual inclusion in this case and its transpose, $\bm{\Pi}^{\top}: V_{\text{DG}}' \mapsto V_{\text{CG}}'$, is the usual $L^2$ projection. 

Now, we can introduce the auxiliary space preconditioner, $\bm{B}_{\text{DG}}: V_{\text{DG}}' \mapsto V_{\text{DG}}$, as follows,
\begin{equation}\label{eqn:ASP-exact}
	\bm{B}_{\text{DG}} = \bm{S}_{\text{DG}} + \bm{\Pi} \bm{A}_{\text{CG}}^{-1} \bm{\Pi}^\top,
\end{equation}
where the operator $\bm{S}_{\text{DG}}: V_{\text{DG}}' \mapsto V_{\text{DG}}$ is the so-called smoother operator, which usually handles the oscillatory high frequency components in $V_{\text{DG}}$.  In this work, we simply use the Jacobi smoother.  Other smoothers, such as the Gauss-Seidel method, could be used here as well. 

As we can see, to implement the preconditioner $\bm{B}_{\text{DG}}$~\eqref{eqn:ASP-exact} in practice, we need to invert $\bm{A}_{\text{CG}}$ exactly.  However, this might be challenging or even impractical when the problem size is large, polynomal degree is high, and/or the anisotropicity is strong.  Therefore, in practice, we usually replace $\bm{A}_{\text{CG}}^{-1}$ by a robust preconditioner $\bm{B}_{\text{CG}}: V_{\text{CG}}' \mapsto V_{\text{CG}}$ and $\bm{B}_{\text{CG}} \approx \bm{A}_{\text{CG}}^{-1}$, i.e., we solve the auxiliary problem in $V_{\text{CG}}$ approximately.  This results in the following inexact version auxiliary space preconditioner,
\begin{equation*}
	\bm{B}^{\text{inexact}}_{\text{DG}} = \bm{S}_{\text{DG}} + \bm{\Pi} \bm{B}_{\text{CG}} \bm{\Pi}^\top.
\end{equation*}

Next we discuss our choices of $\bm{B}_{\text{CG}}$. Notice that $\bm{A}_{\text{CG}}$ is the linear system obtained by using $H^1$-FEM for solving the anisotropic diffusion equation~\eqref{eqn: model-1}-\eqref{eqn:model-2}. When $k=1$, i.e., the linear FEM method, it is well-known that tailored multigrid (MG) methods and their algebriac variants, algebraic multigrid (AMG) methods, can be used~\cite{MV1994,GO1995,S1998,GHT2009,S2012,BCZ2012,WCXX2012,YXZ2013}. Therefore, we use $\bm{B}_{\text{CG}} = \bm{B}_{\text{MG}}$ when $k=1$. In general, many existing efficient solvers for solving anisotropic diffusion problem discretized using linear elements can be applied here. 

For high-order elements, i.e. $k \geq 2$, we again use auxiliary space preconditioning framework to develop an efficient solver by considering the following product of auxiliary spaces for $V_{\text{CG}}$,
\begin{equation*}
	\bar{V}_{\text{CG}} = V_{\text{CG}} \times V_{\text{Linear}},
\end{equation*}
where $V_{\text{Linear}}$ is the linear finite element space, i.e., $V_{\text{CG}}$ with $k=1$. Since $V_{\text{Linear}} \subset V_{\text{CG}}$, this again can be considered as a subspace decomposition $V_{\text{CG}} = V_{\text{CG}} + V_{\text{Linear}}$.  The transfer operator here is just the inclusion $\bm{I}: V_{\text{Linear}} \mapsto V_{\text{CG}}$ and its transpose is the usual $L^2$ projection $\bm{Q} = \bm{I}^{\top}: V_{\text{CG}}' \mapsto V_{\text{Linear}}'$.  To solve the auxiliary problem in $V_{\text{Linear}}$, i.e., $V_{\text{CG}}$ with $k=1$, as we discussed before, we simply apply MG methods $\bm{B}_{\text{MG}}$.  Therefore, the overall auxiliary space preconditioner $\bm{B}_{\text{CG}}$ can be defined as $\bm{B}_{\text{CG}} = \bm{S}_{\text{CG}} + \bm{I} \bm{B}_{\text{MG}} \bm{Q}$ where $\bm{S}_{\text{CG}}$ denotes some smoother in $V_{\text{CG}}$ which we will discus later. Overall, the choice of $\bm{B}_{\text{CG}}$ in our implementation is
\begin{equation*}
	\bm{B}_{\text{CG}} = 
	\begin{cases}
		\bm{B}_{\text{MG}}, \quad  &\text{if} \ k=1, \\
		\bm{S}_{\text{CG}} + \bm{I} \bm{B}_{\text{MG}} \bm{Q}, \quad &\text{if} \ k \ge 2.
	\end{cases}
\end{equation*}

Since we are dealing with anisotropic problems and high-order elements are used, the choice of the smoother $\bm{S}_{CG}$ is crucial to achieve a good performance. As suggested in~\cite{P1992,P1994,ASVZ2017}, when using high-order to solve isotropic problems, Schwarz-type block smoothers should be used and the blocks should be the blocks corresponding to the vertex patchs.  On the other hand, for anisotropic problems, line smoothers should be applied when the standard coarsening strategy is used in the MG methods~\cite{WCXX2012,YXZ2013,TOS2000,PC2012}.  Therefore, we use Schwarz-type block line smoothers as $\bm{S}_{CG}$ and each block includes all the elements sharing vertices along a line, which is usually orthogonal to the anisotropic direction. 

To summarize, the inexact version auxiliary space preconditioner is 
\begin{equation}\label{eqn:ASP-inexact}
	\bm{B}_{\text{DG}}^{\text{inexact}} = 
	\begin{cases}
		\bm{S}_{\text{DG}} + \bm{\Pi}\bm{B}_{\text{MG}} \bm{\Pi}^\top, \quad & \text{if} \ k=1,\\
		\bm{S}_{\text{DG}} + \bm{\Pi} \bm{S}_{\text{CG}} \bm{\Pi}^\top + \bm{\Pi} \bm{I} \bm{B}_{\text{MG}} \bm{Q} \bm{\Pi}^\top, \quad & \text{if} \ k \geq 2.
	\end{cases}
\end{equation}
We will specify our implementations of $\bm{S}_{\text{DG}}$, $\bm{S}_{\text{CG}}$, and $\bm{B}_{\text{MG}}$ for the test problems in the numerical experiments later.

%

%

\section{Numerical Experiment for Testing the Preconditioners}\label{sect:num-preconditioner}
\begin{figure}[H]
\centering
\begin{tabular}{cc}
\includegraphics[width =.3\textwidth]{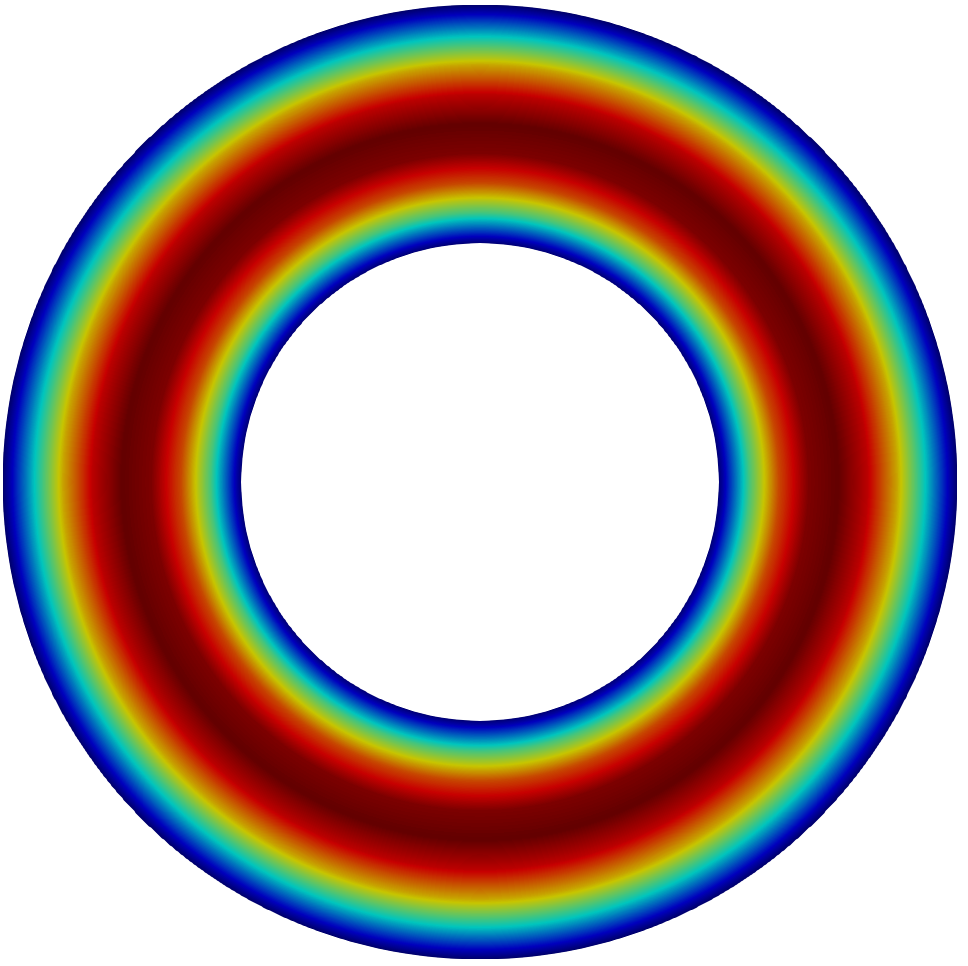}
&\includegraphics[width =.32\textwidth]{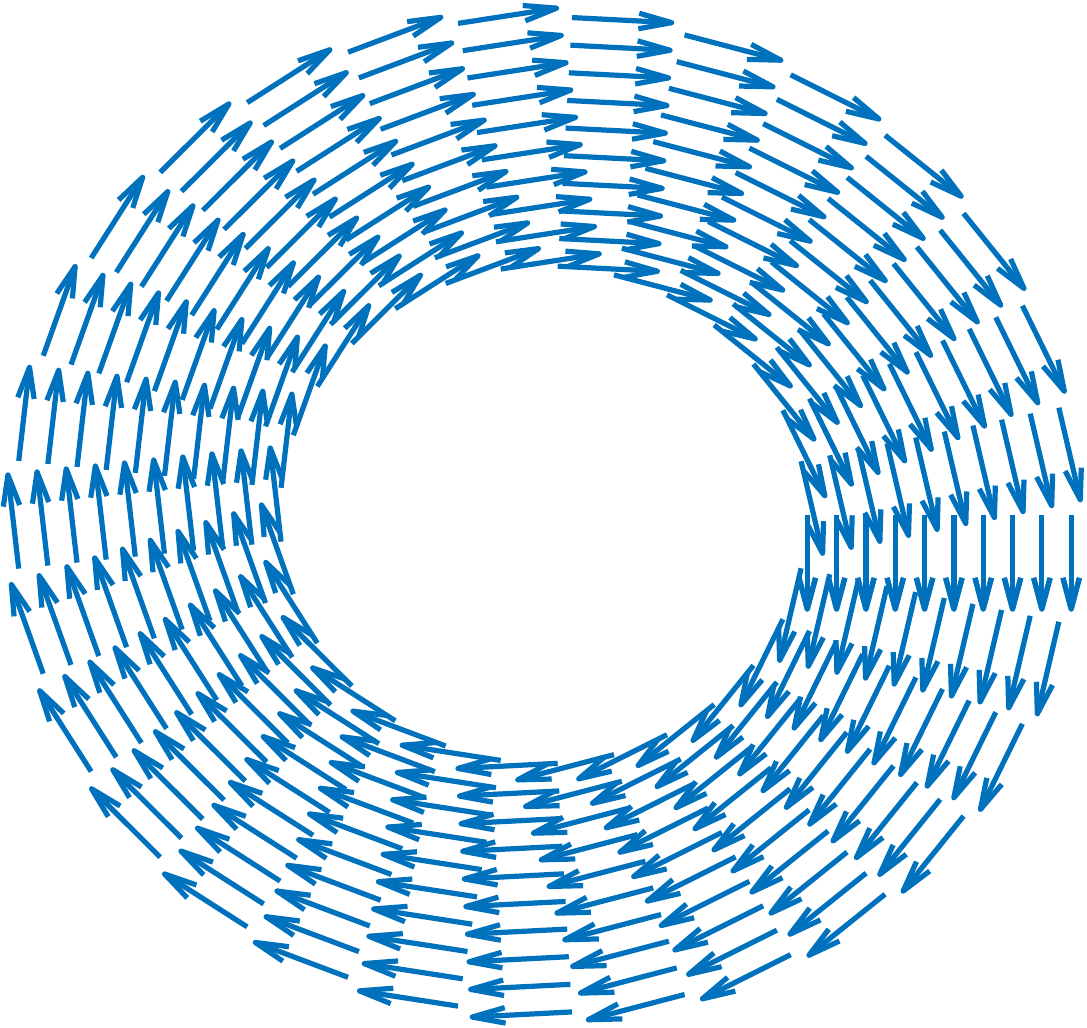}\\
(a) & (b)
\end{tabular}
\caption{Test 5: (a). Plot of Exact Solution. (b) Plot of the magnetic field ${\bf b} = (b_1,b_2)^\top$.}\label{fig:annulus-tri-L2}
\end{figure}

\noindent\textbf{Problem Setting. Test 5.}
In this test, let the computational domain $\Omega$ be an annulus with exterior radius $r = 1$ and interior radius $r = 0.5$. We shall choose the annulus test with following exact solution, diffusion coefficient, and external source as
\begin{eqnarray*}
u &=& \sqrt{\frac{3}{4r}}\sin(2\pi r-\pi), \\
b_1 &=& \frac{y}{r},\ b_2 =- \frac{x}{r}\\
f &=& \sqrt{\frac{3}{4r^5}}(4\pi^2 r^2 - \frac{1}{4})\sin(2\pi r - \pi).
\end{eqnarray*}
The exact solution and the magnetic field ${\bf b}$ are plotted in Figure~\ref{fig:annulus-tri-L2}. The diffusion tensor is of the form (\ref{eq:diff-coef}). We shall test the numerical performance with varying values in $D_{\|}$ and $D_{\perp} = 1$. In the following, we shall demonstrate
\begin{itemize}
\item The high order scheme has an advantage in the both degrees of freedom and conditioning. 
\item The effectiveness of our proposed auxiliary space preconditioners $\bm{B}_{\text{DG}}$ ~\eqref{eqn:ASP-exact} and $\bm{B}_{\text{DG}}^{\text{inexact}}$~\eqref{eqn:ASP-inexact}. 
\end{itemize}

\subsection{Examination of conditioning of the high order schemes for fixed accuracy}

\begin{figure}[H]
\centering
\begin{tabular}{ccc}
\includegraphics[width =.3\textwidth]{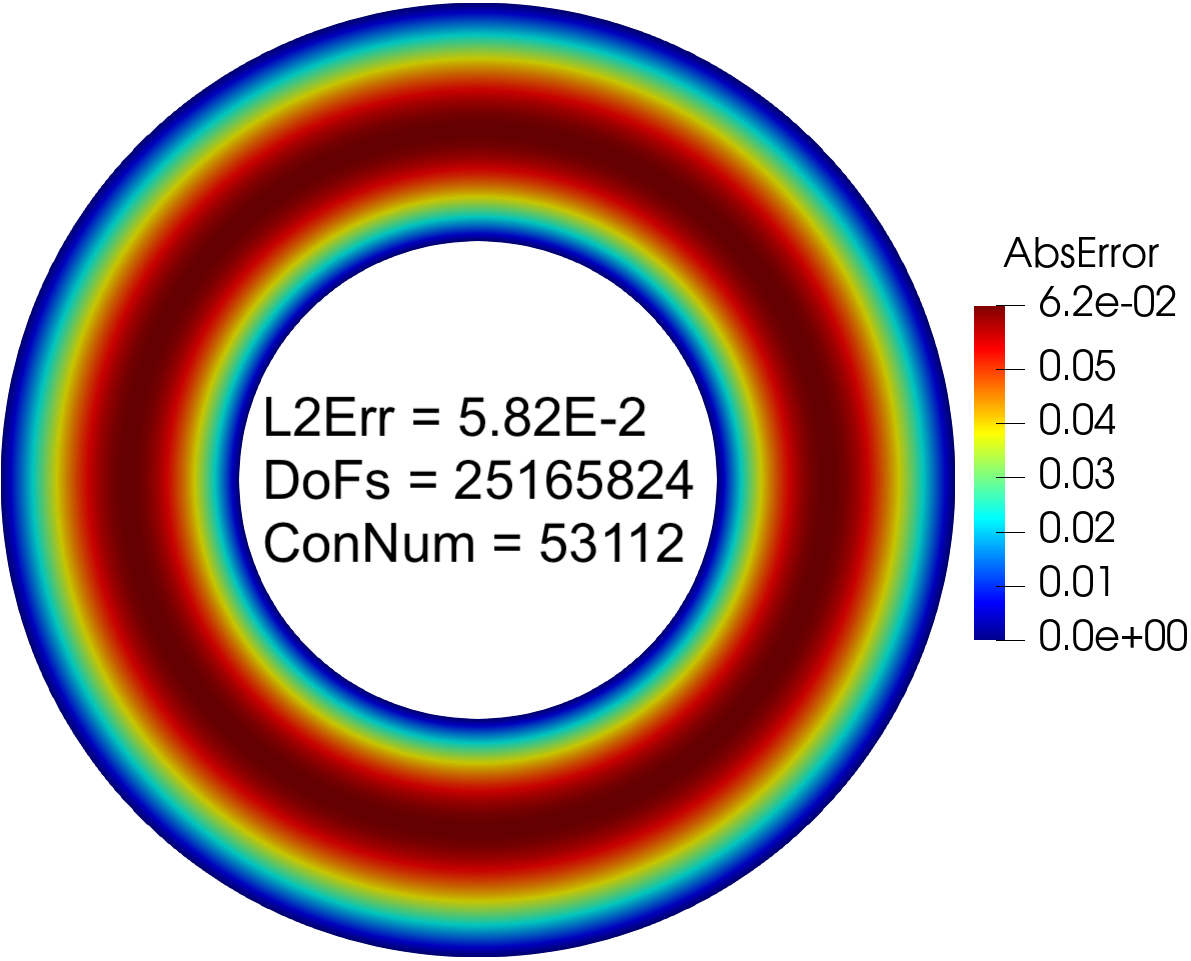}
&\includegraphics[width =.3\textwidth]{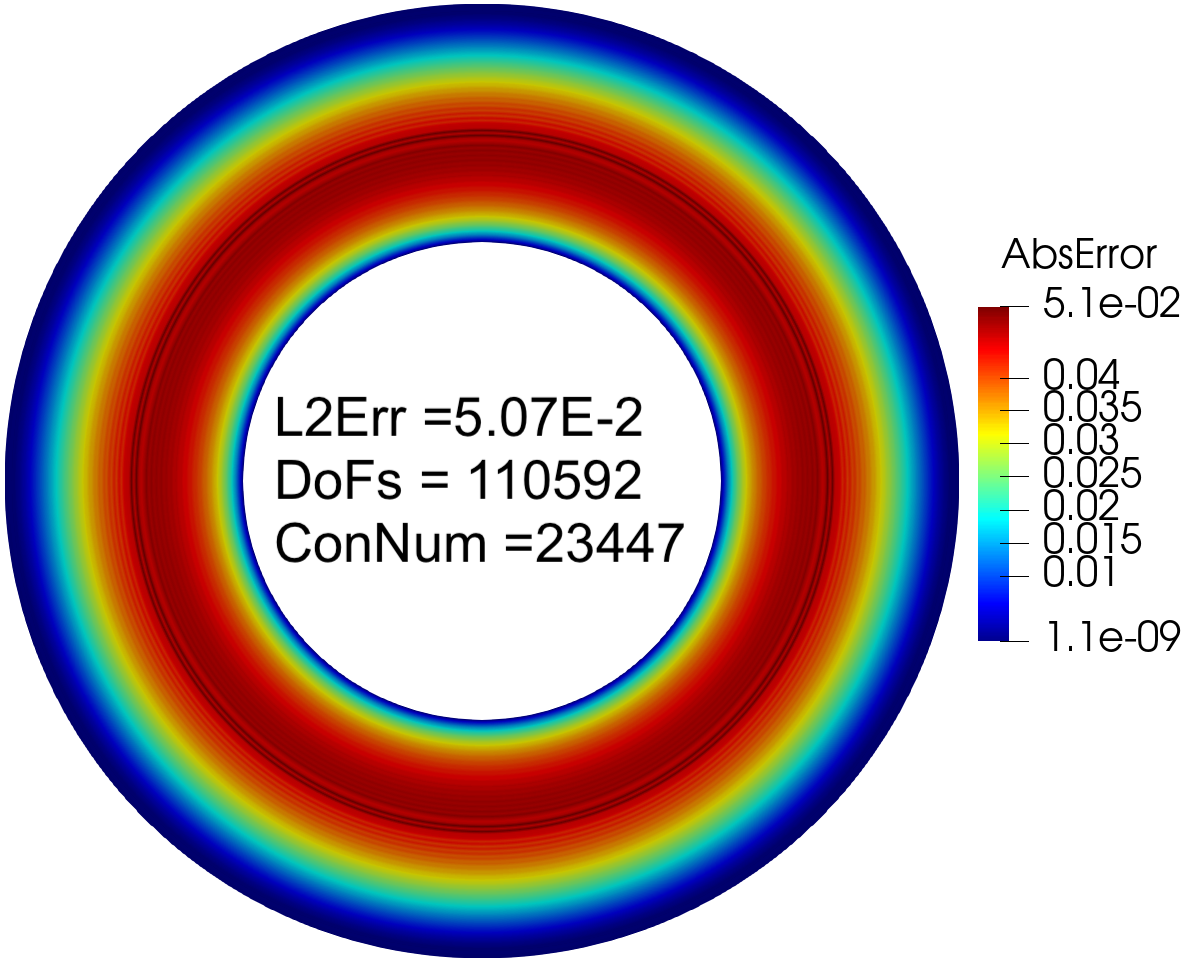}
&\includegraphics[width =.3\textwidth]{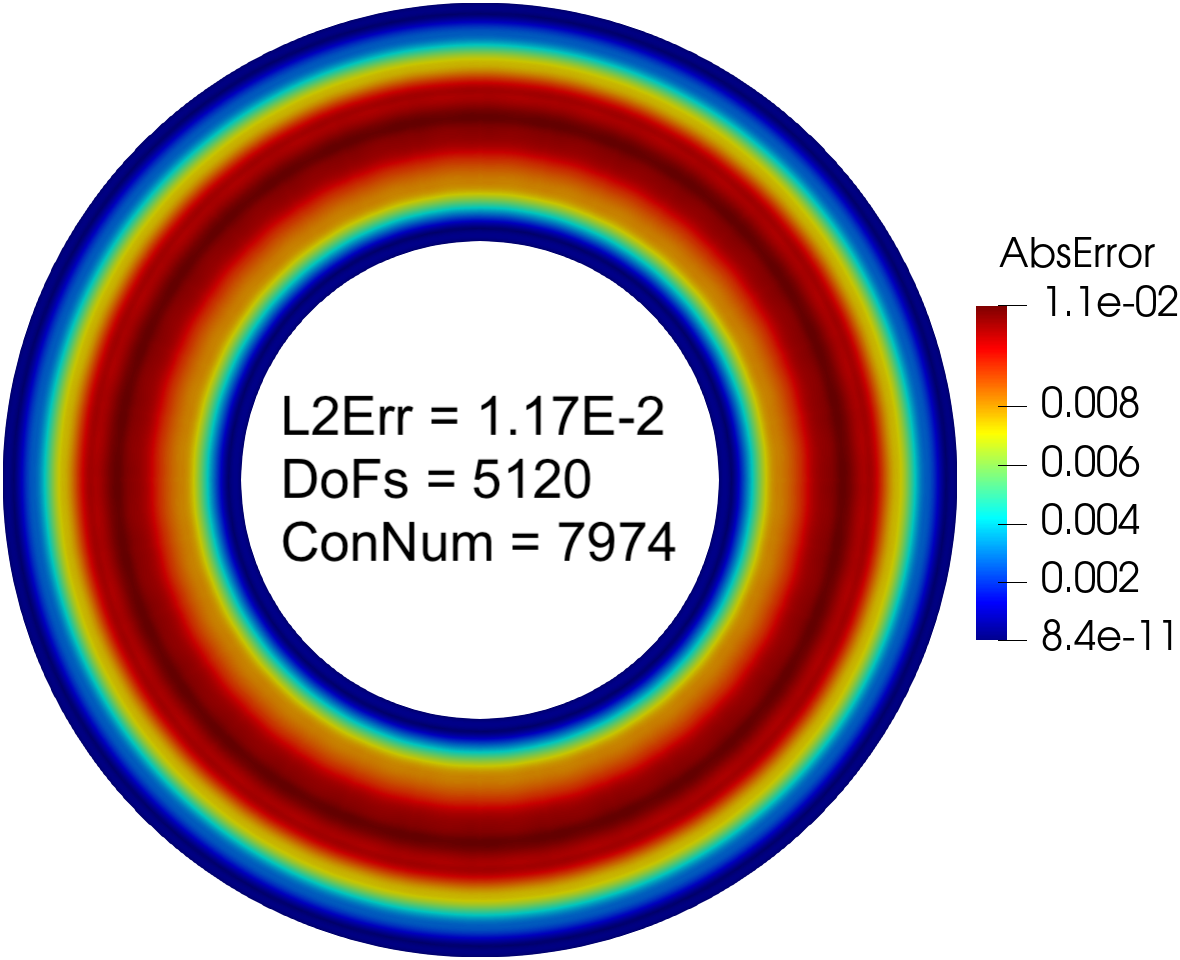}
\\
(a) & (b) &(c)
\end{tabular}
\caption{Plot of errors $|u-u_h|$ for $D_{\|}$ = 1E6 on: (a) $N_r = 1024, k = 1$; (b) $N_r = 48, k = 2$; (c) $N_r = 8, k = 3$.}\label{fig:abserr}
\end{figure}

While in Section~\ref{Num-ConTest} it was seen that the condition number of the resulting linear system increases with scheme order, that test did not hold the accuracy constant and so the as the order was increased, so was the accuracy. In this section we examine the  condition number (ConNum) while holding the accuracy as measured by $L^2$-error (L2Err) approximately constant, and additionally compare the required degrees of freedom (DoFs). Figure~\ref{fig:abserr} plots the absolute errors on different meshes with varying polynomial degrees for the anisotropic ratio $\frac{D_{\|}}{D_\perp} = $1E6. We observe that for linear DG elements, and in order to obtain the $L^2$-error$ = 5.8$E-2, a very fine mesh with $N_r = 1024$ and $N_{\theta} = 4096$ is required. The DoFs of the corresponding linear system are around $2.5$E7. By using quadratic polynomials, we only need $N_r = 48$ and $N_{\theta} = 192$, which requires DoFs$ = 1.1$E5 to achieve a similar $L^2$-error. As we increase the degree to $k = 3$, there is a $21\times$ reductions in DoFs to achieve similar accuracy. Examining the condition numbers, we observe a significant reduction from 5.3E4 for linear elements, to 2.3E4 for quadratic elements, and to 8.0E3 for cubic elements.
This suggests that employing high-order methods not only obtains the required accuracy, but also achieves significant computational savings. 

\subsection{Performance of the preconditioner}
In this section we perform three tests to demonstrate the performance of our proposed preconditioner. In Sec.~\ref{sec:cond-linear-exact} we use different meshes with varying values in anisotropy to examine the preconditioner performance using $\bm{B}_{\text{DG}}$ applied to the linear finite element problem and solving $H^1$ system exactly; in Sec.~\ref{sec:high-order-exact} we look at performance as a function of scheme order, but still using $\bm{B}_{\text{DG}}$ to give an exact solution to the $H^1$ system; and finally in Sec.~\ref{sec:high-order-approx} we present preconditioner performance for varying order while using $\bm{B}_{\text{DG}}^{\text{inexact}}$ to solve the $H^1$ system approximately. In all three cases, our proposed preconditioner is effective. 

\begin{table}[H]
\centering
\begin{tabular}{c|cc|cc}
	\hline \hline 
$D_{\|}$ & Cond($\bm{B}_{\text{DG}}^{-1}\bm{A}_{\text{DG}}$) & PCG Iter($\bm{B}_{\text{DG}}^{-1}\bm{A}_{\text{DG}}$) & Cond($\bm{A}_{\text{DG}}$) & CG Iter($\bm{A}_{\text{DG}}$) \\ \hline \hline 
 &\multicolumn{4}{c}{$N_{r} = 8, N_{\theta} = 32$ }\\  \hline
1 	        &1.93E+01	&19	&6.92E+02	&129\\ 
1E+2	&3.33E+01	&28	&1.95E+04	&596\\
1E+4	&3.28E+01	&32	&4.16E+04	&1411\\
1E+6	&3.35E+01	&33	&4.80E+04	&1453\\
1E+8	&3.30E+01	&33	&4.46E+04	&1375\\
1E+10	&3.37E+01	&36	&4.29E+04	&1289\\ \hline \hline 
&\multicolumn{4}{c}{$N_{r} = 16, N_{\theta} = 64$ }\\ \hline
1 	         &2.16E+01	&17	&2.66E+03	&215\\
1E+2	&3.50E+01	&26	&2.66E+04	&1152\\
1E+4	&3.41E+01	&33	&4.85E+04	&4922\\
1E+6	&3.81E+01	&35	&4.87E+04	&5000+\\
1E+8	&3.79E+01	&35	&5.39E+04	&5000+\\
1E+10	&3.72E+01	&35	&4.74E+04	&5000+\\ \hline \hline 
&\multicolumn{4}{c}{$N_{r} = 32, N_{\theta} = 128$ }\\ \hline
1               &2.28E+01 &25 &7.19E+03 &336 \\
1E+2 & 3.44E+01 &37 &3.66E+04 &1965 \\
1E+4 & 3.76E+01 &40 &5.11E+04 &5000+ \\
1E+6 & 4.11E+01 &39 &5.14E+04 &5000+\\
1E+8 & 4.19E+01 &41 &5.14E+04 &5000+\\
1E+10 &4.40E+01 &42 &5.14E+04 &5000+\\ \hline\hline
\end{tabular}
\caption{Conditioning results for $k = 1$ on triangular mesh}\label{Table:ASMP-k1}
\end{table}



\subsubsection{Conditioning Results for Linear Finite Elements}
\label{sec:cond-linear-exact}
Here we demonstrate the effectiveness of our proposed preconditioner for linear finite elements with varying mesh size spanning $N_r = 8, 16, 32$ with $N_\theta = 32, 64, 128$. The conditioning results are reported in Table~\ref{Table:ASMP-k1}. The condition numbers are again computed by the command ``{\it condest}'' in PyAMG and the conjugate gradient (CG) method with stopping criterion ($\epsilon<$1E-6) has been applied for validating the effectiveness of our preconditioner.  Here, we use Jacobi method for $\bm{S}_{DG}$ and, for the $H^1$-problem, the action of $\bm{A}_{\text{CG}}^{-1}$ in~\eqref{eqn:ASP-exact} is obtained by adopting the direct solver. 
The results show that the condition number for the stiffness matrix $\bm{A}_{\text{DG}}$ is related to the mesh size and the value of $D_{\|}$. When refining the mesh, the condition number usually increases at about $\mathcal{O}(h^{-2})$. However, the results show that the condition number from a fixed mesh with increasing values in $D_{\|}$ is bounded by the order $10^{4}$. We do not know the reason for this behaviour. In contrast to the results for original stiffness matrix, the pre-conditioned system has a bounded condition number, which is around $40$.  When the CG iterative solver is applied to solve the problem, the pre-conditioned system needs around 40 iterations to achieve the required accuracy 1E-6. We observe that the iteration numbers for CG (with tol = 1E-6) for the preconditioned system are almost independent with the values in $D_{\|}$ and $N_r, N_{\theta}$. However, the required iteration numbers in the CG for the original linear system are significantly larger where when the mesh is finer or the value in $D_{\|}$ is large, the CG solver fails to converge within 5000 iterations. This validates the efficiency and effectiveness of our proposed preconditioner.

\subsubsection{Solve $H^1$ problem exactly for varying scheme order: $\bm{B}_{\text{DG}}$.}
\label{sec:high-order-exact}
%
For varying degree of the numerical schemes, here we test the performance of $\bm{B}_{\text{DG}}$~\eqref{eqn:ASP-exact} where the $H^1$-problem, $\bm{A}_{\text{CG}}$ is solved exactly via a direct solver. In addition, we use the Jacobi method for $\bm{S}_{\text{DG}}$.  The numerical results are shown in Table~\ref{tab:exact_H1} and Table~\ref{tab:exact_H1_v2} for two different mesh sizes.  We used the generalized minimal residual method (GMRES) with stopping tolerance 1E-6 to solve the pre-conditioned system in order to achieve more robust performance since in practice the linear system becomes ill-conditioned and CG may suffer loss of orthogonality due to numerical instability. As shown in these two tables, the required number of iterations is almost constant and independent of polynomial degree $k$, values of $D_{\|}$, or the mesh size $1/N_r\times1/N_\theta$. The convergence results for the relative residual in the iterative solver for $k = 2$ are plotted in Figure~\ref{fig:annulus-tri-residual}a and Figure~\ref{fig:annulus-tri-residual_v2}a. The reduction rates in this plot is at the rate $0.8$.
Therefore, our results demonstrate that the $H^1$-problem can provide an effective and robust preconditioner which controls the condition number in the DG system. However, in the cases of finer mesh size, high order schemes, and/or extreme anisotropy, the direct solver for the corresponding $H^1$-problem may still be infeasible. Hence, the inexact version $\bm{B}_{\text{DG}}^{\text{inexact}}$, in which the same $H^1$-problem is solved approximately, is adopted below.

\subsubsection{Solve $H^1$ problem inexactly for varying scheme order: $\bm{B}_{\text{DG}}^{\text{inexact}}$.}
\label{sec:high-order-approx}
\begin{figure}[H]
 \centering
 \begin{tabular}{cc}
 \includegraphics[width =.45\textwidth]{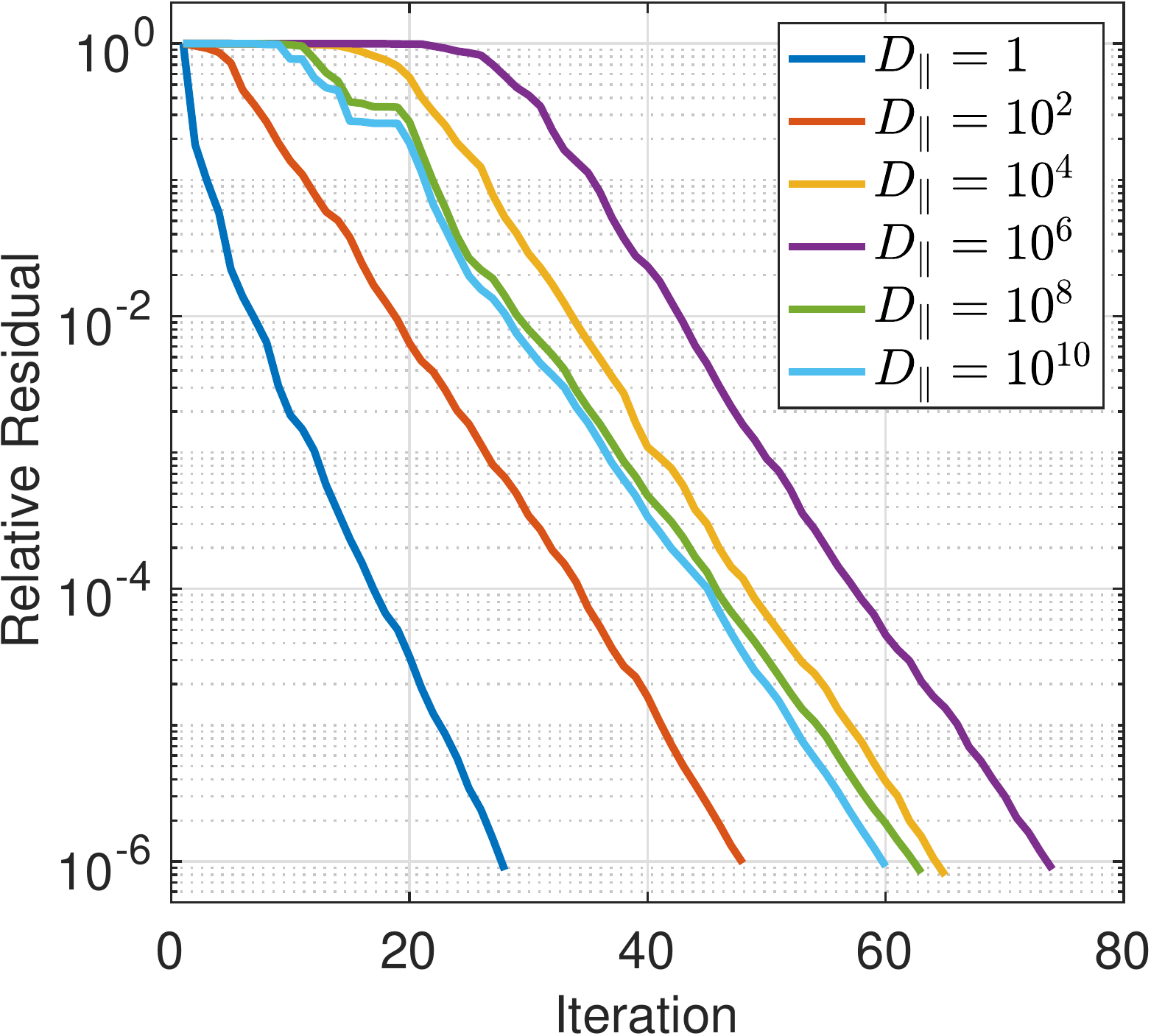} 
 &\includegraphics[width =.45\textwidth]{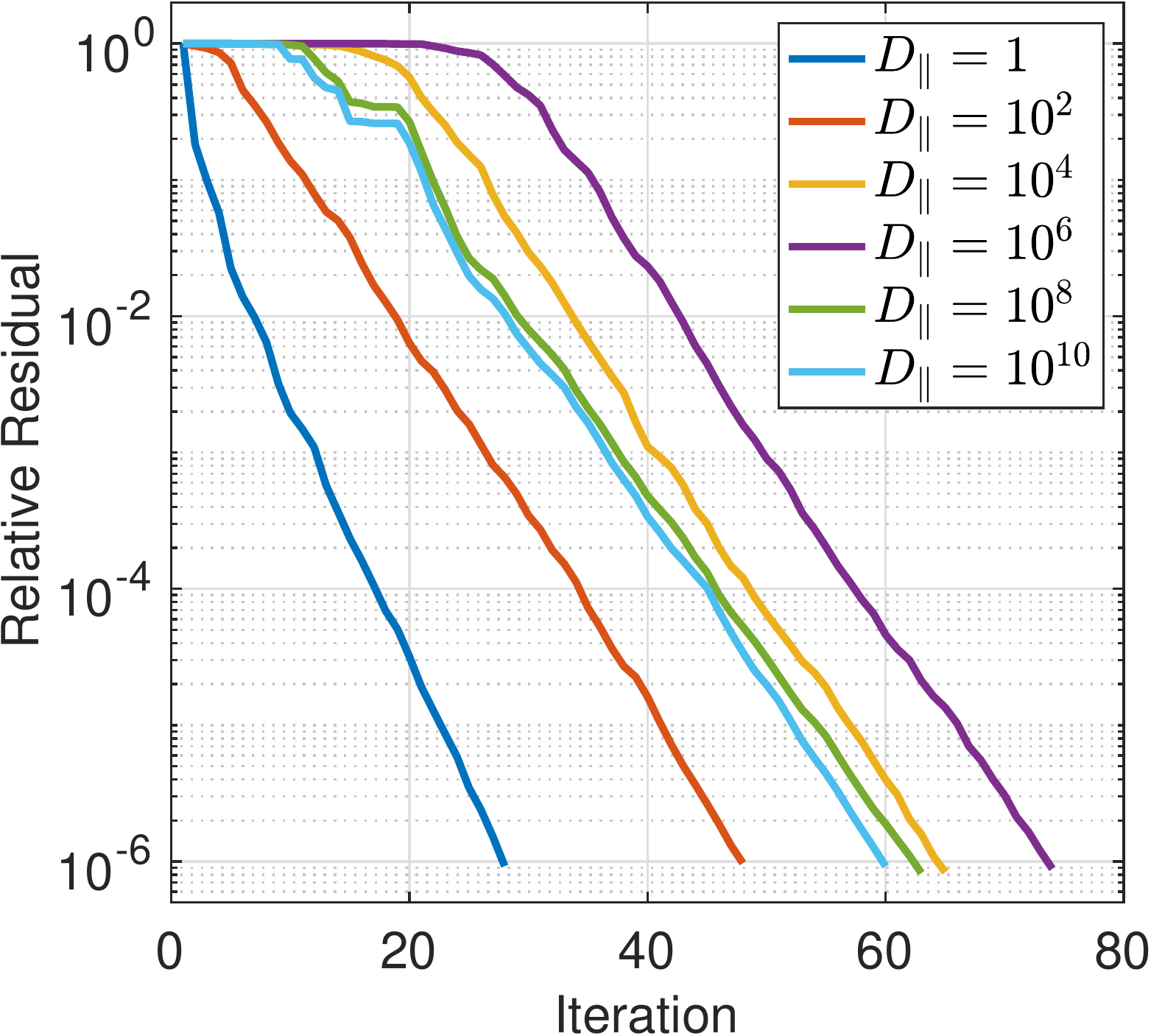} 
 \\
 (a) & (b)
 \end{tabular}
 \caption{Plot of relative residual $\frac{\|\bm{r}\|}{\|\bm{r}_0\|}$ for the iterative solver for $k = 2$ on the mesh with $N_r = 8, N_\theta = 32$: (a) Solve the $H^1$ problem exactly; (b) Solve the $H^1$ problem inexactly.}\label{fig:annulus-tri-residual}
 \end{figure}

\begin{figure}[H]
 \centering
 \begin{tabular}{cc}
 \includegraphics[width =.45\textwidth]{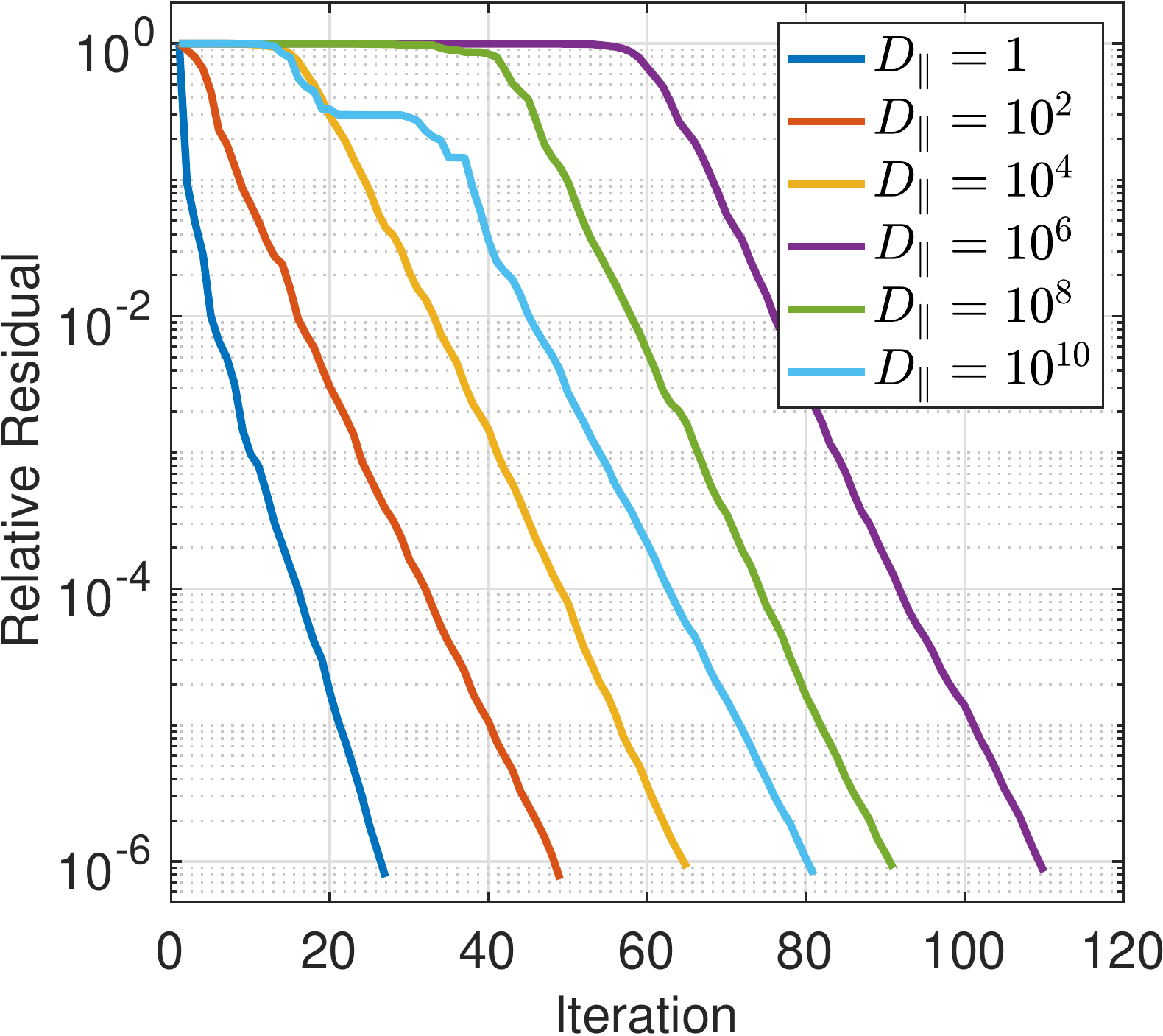} 
 &\includegraphics[width =.45\textwidth]{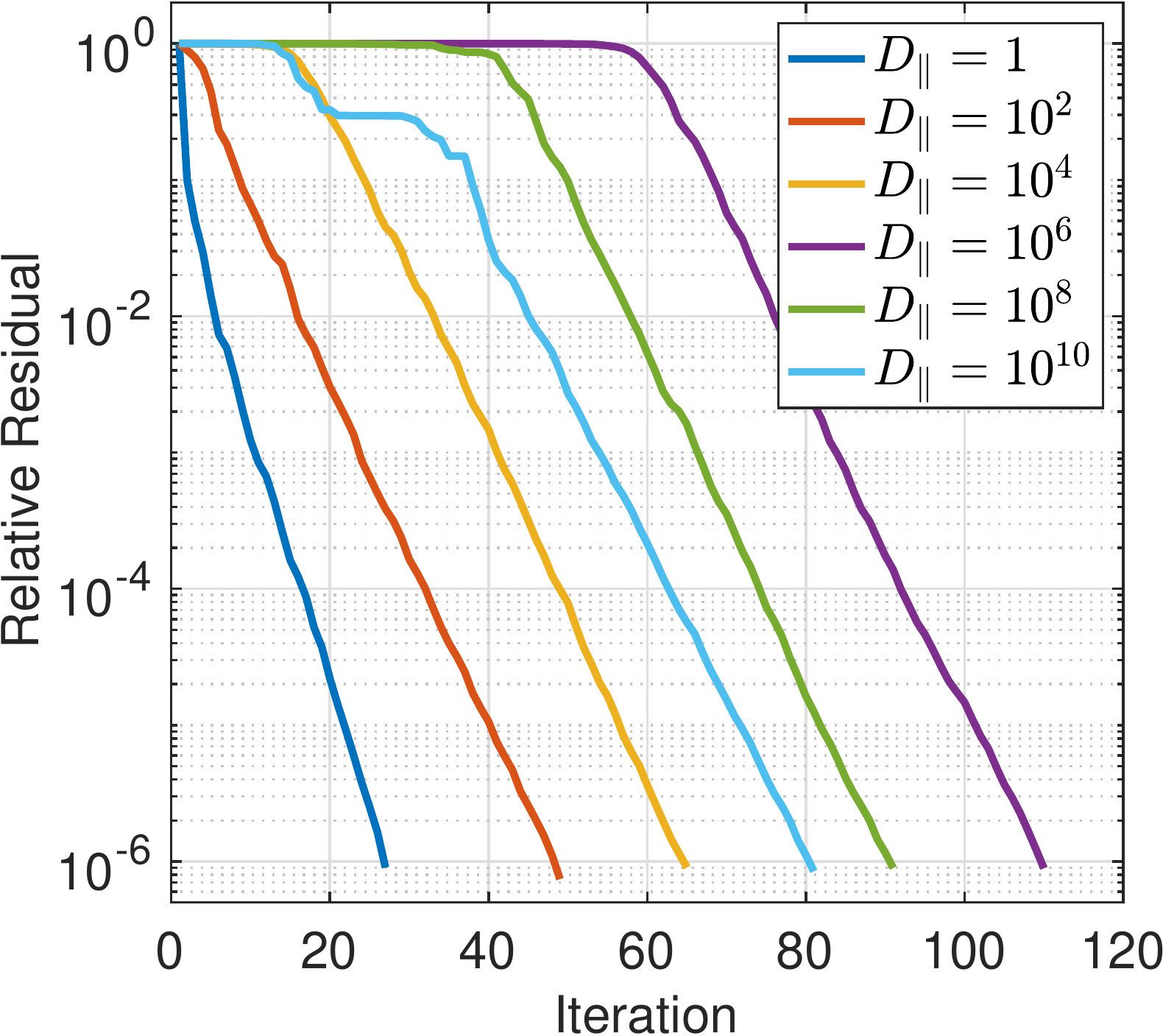} 
 \\
 (a) & (b)
 \end{tabular}
 \caption{Plot of relative residual $\frac{\|\bm{r}\|}{\|\bm{r}_0\|}$ decaying for the iterative solver for $k = 2$ on the mesh with $N_r = 16, N_\theta = 64$: (a) Solve the $H^1$ problem exactly; (b) Solve the $H^1$ problem inexactly.}\label{fig:annulus-tri-residual_v2}
 \end{figure}

		In our implementation of $\bm{B}_{\text{DG}}^{\text{inexact}}$, we again use the Jacobi method for $\bm{S}_{\text{DG}}$. We use the smoothed aggregation AMG method for $\bm{B}_{\text{MG}}$ when solving the auxiliary problem in the linear finite element space as suggested in~\cite{S2012}.  As we discussed, $\bm{S}_{\text{CG}}$ should be a Schwarz-type block line smoother. To be more precise, for this test problem, and since the computational domain is an annulus and the anisotropicity is along the circular direction, the blocks for the line smoother should be constructed along the radial direction. In practice, we first find the vertices along each radial direction on the mesh and then use all the unknowns (which can be found by looking at the nonzeros in the stiffness matrix $\bm{A}_{\text{CG}}$) that are connected to those vertices to form line blocks along the radial direction.  We refer to the block line smoother using those blocks as $\bm{S}_{\text{CG}}^{\text{rad}}$.  In addition, due to the fact the domain is an annulus, our test problem is periodic in the circular direction which leads to extra high-frequency components along the circular direction. Therefore, we also use block line smoothers along the circular direction. The line smoother in the circular direction is more precisely described by the following procedure: we find the vertices along each circular direction on the mesh and use all the unknowns (which, again, can be found by looking at the nonzeros in $\bm{A}_{\text{CG}}$) that are connected to those vertices to form line blocks along the circular direction. We refer to the block line smoother using those blocks as $\bm{S}_{\text{CG}}^{\text{cir}}$.  Overall, we use $\bm{S}_{\text{CG}}^{\text{rad}}$ and $\bm{S}_{\text{CG}}^{\text{cir}}$ in a multiplicative fashion to define $\bm{S}_{\text{CG}}$, i.e., 
		$\bm{S}_{\text{CG}} := \bm{S}_{\text{CG}}^{\text{cir}} + \bm{S}_{\text{CG}}^{\text{rad}} - \bm{S}_{\text{CG}}^{\text{cir}}\bm{A}_{\text{CG}}\bm{S}_{\text{CG}}^{\text{rad}}$.  In addition, to further improve the robustness of the preconditioner, we use $\bm{B}_{\text{CG}}$ as a preconditioner for GMRES to solve $\bm{A}_{\text{CG}}$ inexactly with a tolerance of $1E-2$ and switch the outer Krylov method for solving $\bm{A}_{\text{DG}}$ to a flexible version of GMRES with a tolerance of $1E-6$. 
		 The numerical results are presented in Table~\ref{tab:approx_H1} and~\ref{tab:approx_H1_v2} for mesh size of $N_r = 8, \ N_{\theta}=32$ and $N=16, \ N_{\theta}=64$, respectively. As we can see, the number of iterations roughly stay the same for the different polynomials degrees $k$ and anisotropic diffusion coefficient $D_{\|}$. In addition, since we solve the auxiliary $H^1$ problem inexactly, the Krylov method takes slightly more iterations to converge than when the $H^1$ problem is solved exactly, which is expected. The convergence results for the relative residual of the iterative solve for $k = 2$ are plotted in Figure~\ref{fig:annulus-tri-residual}b and Figure~\ref{fig:annulus-tri-residual_v2}b.  Overall, the numerical results demonstrate that the inexact preconditioner is robust and effective for this test problem. 
		


\begin{table}[H]
    \centering
    \begin{tabular}{c||c|c|c|c|c|c|c|c}
    \hline \hline 
      $D_{\|}$  & $k=1$ & $k=2$ & $k=3$ & $k=4$ & $k=5$ & $k=6$ & $k=7$ & $k=8$ \\ \hline \hline
       1          &  23  & 27  & 22  & 24 & 23 &28   &26  &30  \\
       1E+2   &  39  & 47  & 39  & 40 & 34 &38   &37    &42 \\
       1E+4   &  47  & 64  & 52  & 55 & 51 &53   &54    &56 \\
       1E+6   &  46  & 73  & 56  & 58 & 53 &54   &55    &56 \\
       1E+8   &  39  & 62  & 59  & 58 & 53 &54   &55    &56 \\
       1E+10  & 35      &59     & 52  & 54 &53  &54  &55    &56 \\
       \hline \hline 
    \end{tabular}
    \caption{Performance of~$\bm{B}_{\text{DG}}$ (solve the $H^1$ problem exactly) for $N_r = 8$ and $N_{\theta} = 32$}
    \label{tab:exact_H1}
\end{table}


\begin{table}[H]
    \centering
    \begin{tabular}{c||c|c|c|c|c|c|c|c}
    \hline \hline 
      $D_{\|}$  & $k=1$ & $k=2$ & $k=3$ & $k=4$ & $k=5$ & $k=6$ & $k=7$ & $k=8$ \\ \hline \hline
       1          &  25  & 27  & 23  & 25 & 27  &30  &32  &35  \\
       1E+2   &  39  & 47  & 39  & 40 &38  &43   &46    &49 \\
       1E+4   &  47  & 64  & 53  & 56 &53  &57  &60    &62 \\
       1E+6   &  46  & 73  & 57  & 58 &54  &54  &55    &57 \\
       1E+8   &  39  & 62  & 63  & 58 &53  &54  &55    &56 \\
       1E+10  &35      &59     & 56  & 60 &60  &54  &55    &56 \\
       \hline \hline 
    \end{tabular}
    \caption{Performance of~$\bm{B}_{\text{DG}}^{\text{inexact}}$ (solve the $H^1$ problem inexactly) for $N_r = 8$ and $N_{\theta} = 32$}
    \label{tab:approx_H1}
\end{table}

\begin{table}[H]
    \centering
    \begin{tabular}{c||c|c|c|c|c|c|c|c}
    \hline \hline 
      $D_{\|}$  & $k=1$ & $k=2$ & $k=3$ & $k=4$ & $k=5$ & $k=6$ & $k=7$ & $k=8$ \\ \hline \hline
       1          & 22   &26     &20   &23  &21  &25   &22    &27  \\
       1E+2   &40    &48     &35   &35  &27  &36   &29    &37 \\
       1E+4   &52    &64     &49   &51  &42  &45   &47    &50 \\
       1E+6   &54    &109   &60   &62  &53  &55   &55    &57 \\
       1E+8   &53    &90     &60   &62  &53  &55   &55    &56 \\
       1E+10  &43   &80     &75   &62  &53  &55   &55    &56 \\
       \hline \hline 
    \end{tabular}
    \caption{Performance of~$\bm{B}_{\text{DG}}$ (solve the $H^1$ problem exactly) for $N_r = 16,N_\theta = 64$}
    \label{tab:exact_H1_v2}
\end{table}

\begin{table}[H]
    \centering
    \begin{tabular}{c||c|c|c|c|c|c|c|c}
    \hline \hline 
      $D_{\|}$  & $k=1$ & $k=2$ & $k=3$ & $k=4$ & $k=5$ & $k=6$ & $k=7$ & $k=8$ \\ \hline \hline
       1          &  22  &26     &23     &25  &28  &31   &33  &35  \\
       1E+2   &40    &48     &35     &36  &35  &40   &41    &45 \\
       1E+4   &52    &64     &50     &52  &50  &61   &65    &69 \\
       1E+6   &54    &109   &60     &63  &61  &58   &58    &62 \\
       1E+8   &53    &90     &65     &63  &57  &60   &64    &64 \\
       1E+10  &43   &80     &80   &80  &61  &69  &69    &69 \\
       \hline \hline 
    \end{tabular}
    \caption{Performance of~$\bm{B}_{\text{DG}}^{\text{inexact}}$ (solve the $H^1$ problem inexactly) for $N_r = 16,N_\theta = 64$}
    \label{tab:approx_H1_v2}
\end{table}

%

\section{Conclusions}\label{sect:conclusion}
In this paper, an interior penalty discontinuous Galerkin finite element scheme has been presented for the discretization of diffusion equations with strong anisotropy. With a high order scheme, the method is potentially able to discretize problems of relevance to magnetized plasma physics where complex magnetic topologies and plasma-facing component geometries (which describe the domain boundary) make the generation of magnetic field aligned meshes difficult or impossible. The preconditioner, constructed by the $H^1$-problem together with the Jacobi smoother, has been demonstrated to solve the corresponding linear systems efficiently. Future applications of our proposed algorithm will include more complicated geometries, including three dimensional domains and Tokamak device shapes.
%
%

\bibliographystyle{elsarticle-num}

\end{document}